\documentclass[12pt]{amsart}

\usepackage{graphicx}

\usepackage{diagrams}

\usepackage{xcolor}

\usepackage{color}

\usepackage{hyperref}
\hypersetup{colorlinks=false,linkbordercolor=red,linkcolor=green,pdfborderstyle={/S/U/W 1}}

\newcommand{\lab}{\label}

\numberwithin{equation}{section}

\newcommand{\ben}{\begin{enumerate}}
\newcommand{\een}{\end{enumerate}}

\newcommand{\bea}{\begin{eqnarray}}
\newcommand{\ba}{\begin{array}}
\newcommand{\bean}{\begin{eqnarray*}}
\newcommand{\ea}{\end{array}}
\newcommand{\eea}{\end{eqnarray}}
\newcommand{\eean}{\end{eqnarray*}}
\newcommand{\beq}{\begin{equation}}
\newcommand{\eeq}{\end{equation}}
\newcommand{\bthm}{\begin{thm}}
\newcommand{\ethm}{\end{thm}}
\newcommand{\blem}{\begin{lem}}
\newcommand{\elem}{\end{lem}}
\newcommand{\bprop}{\begin{prop}}
\newcommand{\eprop}{\end{prop}}
\newcommand{\bcor}{\begin{cor}}
\newcommand{\ecor}{\end{cor}}
\newcommand{\bdfn}{\begin{dfn}}
\newcommand{\edfn}{\end{dfn}}
\newcommand{\brem}{\begin{rem}}
\newcommand{\erem}{\end{rem}}
\newcommand{\bpf}{\begin{proof}}
\newcommand{\epf}{\end{proof}}
\newcommand{\bfact}{\begin{fact}}
\newcommand{\efact}{\end{fact}}
\newcommand{\bobs}{\begin{obs}}
\newcommand{\eobs}{\end{obs}}
\newcommand{\bexam}{\begin{exam}}
\newcommand{\eexam}{\end{exam}}
\newcommand{\bclaim}{\begin{claim}}
\newcommand{\eclaim}{\end{claim}}

\newtheorem{thm}{Theorem}[section]
\newtheorem{prop}[thm]{Proposition}
\newtheorem{lem}[thm]{Lemma}

\newtheorem{cor}[thm]{Corollary}
\newtheorem{dfn}[thm]{Definition}
\newtheorem{rem}[thm]{Remark}
\newtheorem{fact}[thm]{Fact}
\newtheorem{claim}[thm]{Claim}
\newtheorem{obs}[thm]{Observation}
\newtheorem{exam}[thm]{Example}

\newtheorem*{condition'}{Condition 2'}

 \newtheoremstyle{claimstyle}%
   {}
   {}
   {\normalfont}
   {}
   {\itshape}
   {.}
   { }
   {\thmnote{#3}}

\theoremstyle{claimstyle}

\alph{enumii} \roman{enumiii}

\def\cA{\mathcal A}             \def\cB{\mathcal B}       \def\cC{\mathcal C}
\def\cH{\mathcal H}             \def\cF{\mathcal F}       \def\cG{\mathcal G}
           \def\cM{\mathcal M}        \def\cP{{\mathcal P}}
             \def\cV{\mathcal V}       \def\cJ{\mathcal J}   \newcommand{\J}{\mathcal{J}}
\def\cR{\mathcal R}              \def\cS{\mathcal S}

 \def\cD{\mathcal D}   

                \def\Z{{\mathbb Z}}      \def\R{{\mathbb R}}
\def\C{{\mathbb C}}                      \def\oc{{\hat \C}}

\def\F{{\mathcal F}}

\def\a{\alpha}                \def\b{\beta}             \def\d{\delta}
\def\De{\Delta}                          
\def\g{\gamma}                \def\Ga{\Gamma}           \def\l{\lambda} 
                         \def\Om{\Omega}
               \def\sg{\sigma}
\def\Th{\Theta}                          
\def\ka{\kappa}

\newcommand{\ep}{\varepsilon}
\newcommand{\ph}{\varphi}
\newcommand{\al}{\alpha}
\newcommand{\ga}{\gamma}

\def\1{1\!\!1}

\def\and{\text{ and }}

        \def\diam{\text{\rm {diam}}}
\def\dist{\text{{\rm dist}}}

\def\H{\text{{\rm H}}}     \def\HD{\text{{\rm HD}}}

\def\Int{\text{{\rm Int}}}
         \def\P{\text{{\rm P}}}     \def\Id{\text{{\rm Id}}}
 
\def\Hdim{\text{{\rm Hdim}}} \def\HypDim{\text{{\rm HypDim}}}

\def\({\bigl(}                \def\){\bigr)}
\def\lt{\left}                \def\rt{\right}

                        \def\^{\tilde}

            \def\sms{\setminus}
\def\sbt{\subset}             \def\spt{\supset}
\def\gek{\succeq}             \def\lek{\preceq}

\def\comp{\asymp}
           \def\downto{\searrow}
\def\sp{\medskip}             \def\fr{\noindent}        
\def\ov{\overline}            

\def\fr{\noindent}

\def\D{{\mathbb D}}

\def\${$ \displaystyle }

\newcommand{\pft}{{\mathcal{L}}_t}

\newcommand{\Tract}{\Om}

\def\den{\rho}

\newcommand{\pf}{{\mathcal{L}}}

\newcommand{\npf}{\mathcal{\hat L}}

\newcommand{\jul}{\mathcal J}

\def\rad{\J_r}

\def\tstar{\Theta_f}
\def\tstarF{\Theta_F}
\def\tstarj{\Theta_{f,j}}

\begin{document}

\title[Thermodynamic formalism and integral means spectrum]{Thermodynamic formalism and \\ integral means spectrum of logarithmic tracts \\ for transcendental entire functions}


\author[\sc Volker Mayer]{\sc Volker Mayer}
\address{Volker Mayer,
Universit\'e de Lille, UFR de
  Math\'ematiques, UMR 8524 du CNRS, 59655 Villeneuve d'Ascq Cedex,
  France} \email{volker.mayer@univ-lille.fr \newline
  \hspace*{0.42cm} \it Web: \rm math.univ-lille1.fr/$\sim$mayer}

\author[\sc Mariusz Urba\'nski]{\sc Mariusz Urba\'nski}
\address{Mariusz Urba\'nski, 
Department of
  Mathematics, University of North Texas, Denton, TX 76203-1430, USA}
\email{urbanski@unt.edu \newline \hspace*{0.42cm} \it Web: \rm
  www.math.unt.edu/$\sim$urbanski}

\date{\today} \subjclass{111}


\begin{abstract} 
We provide an entirely new approach to the theory of thermodynamic formalism for  entire functions of bounded type.
The key point is that we introduce an integral means spectrum for logarithmic tracts which takes care of the fractal behavior of the boundary of the tract
near infinity. It turns out that this spectrum behaves well as soon as the tracts have some sufficiently nice geometry which, for example, is the case for quasidisk, John or H\"older tracts. In these cases we get a good control of the corresponding transfer operators, leading to full thermodynamic formalism along with its applications such as 
exponential decay of correlations, central limit theorem and a Bowen's formula for the Hausdorff dimension of radial Julia sets.

This approach covers all entire functions for which thermodynamic formalism has been so far established and goes far beyond. It applies in particular to every hyperbolic function from any Eremenko-Lyubich analytic family of Speiser class $\cS$ provided this family contains at least one function with H\"older tracts. The latter is, for example,  the case if the family contains a Poincar\'e linearizer.
\end{abstract}

\maketitle

\section{Introduction} 
The dynamics of a holomorphic function heavily depends on the behavior of the singular set.
The singular set $S(f)$ of an entire function $f:\C\to\C$ is the closure of the set of critical values and finite asymptotic values of $f$.
Eremenko-Lyubich \cite{EL92} introduced and studied class $\cB$ consisting of all entire functions with
bounded singular sets. It has as a subclass Speiser class $\cS$ consisting of entire functions with finite singular sets. In this paper we develop the full theory of thermodynamic formalism for a large collection of entire functions in Eremenko-Lyubich class $\cB$.

When developing the thermodynamic formalism for transcendental functions, one encounters immediately two 
major difficulties: one has to deal with the essential singularity at infinity and to check whether the transfer operator,
which is given by an infinite series, is well defined, i.e. converges, and has sufficiently good properties. 

The first work on thermodynamic formalism
for transcendental functions is due to  Bara\'nski \cite{Baranski95} who considered the tangent family.
Other specific, mainly periodic, functions have been treated in the sequel, see for example \cite{KU02}, \cite{UZ03} and \cite{UZ04}.
The first and, up to now, the only unified approach appeared in \cite{MyUrb08} and in \cite{MUmemoirs}. 
These papers deal with a large class of functions that satisfy a condition on the derivative called the balanced growth condition.
The key point there was to employ Nevanlinna Theory and to make a judicious choice of Riemannian metric.
Here we keep this choice of metric but then we proceed totally differently 
avoiding any use of Nevanlinna Theory. By introducing the integral means spectrum for logarithmic tracts we built the
theory of thermodynamic formalism for many other entire functions from class $\cB$.

\medskip

The main object of this paper is to show that the transfer operator behaves well depending on the geometry of the logarithmic tracts over infinity. 
Consider $f\in \cB$ and suppose that the bounded set $S(f)$ is contained in the unit disk. Then, the components $\Om_j$ of $f^{-1} (\{|z|>1\})$ are the tracts, in fact logarithmic tracts of $f$ over infinity. We assume that there are only finitely many of them: see the definition of class $\cD$ in the next section. 

Let us consider here in this introduction the case where $f^{-1} (\{|z|>1\})$ consists of only one tract $\Om$. Then, $f_{|\Om}$ has the particular form 
$
f=e^\tau,
$
with $\tau $ a conformal map from $\Om$ onto the half plane $\cH=\{\Re z>0\}$  such that 
\beq\label{2_2017_06_28}
\ph:=\tau^{-1}:\cH\to\Om
\eeq
extends continuously to infinity \cite{EL92}.
Although $\partial \Om$ is an analytic curve, near infinity it often resembles more and more a fractal curve.
Typically, going to infinity on $\partial \Om$ is like considering Green lines that are closer and closer to the boundary of possibly fractal domains. To make this precise, consider the rectangles
\beq \label{2}
Q_T:=\big\{\xi\in\C:0<\Re \xi <4T \  \  {\rm and } \  \  -4T < \Im \xi < 4T\big\} 
\eeq
and then the domains
$$
\Omega _T:= \ph ( Q_T )\; , \  \  T\geq 1\,.
$$
The domains $\Omega_T$ form natural exhaustions of $\Om$ and the fractality near infinity of $\partial \Om$ can be observed by considering $\Omega_T$ rescaled by the factors $1/|\ph (T)|$ as $T\to \infty$.

\begin{figure}[h]
   \includegraphics[width=12.5cm]{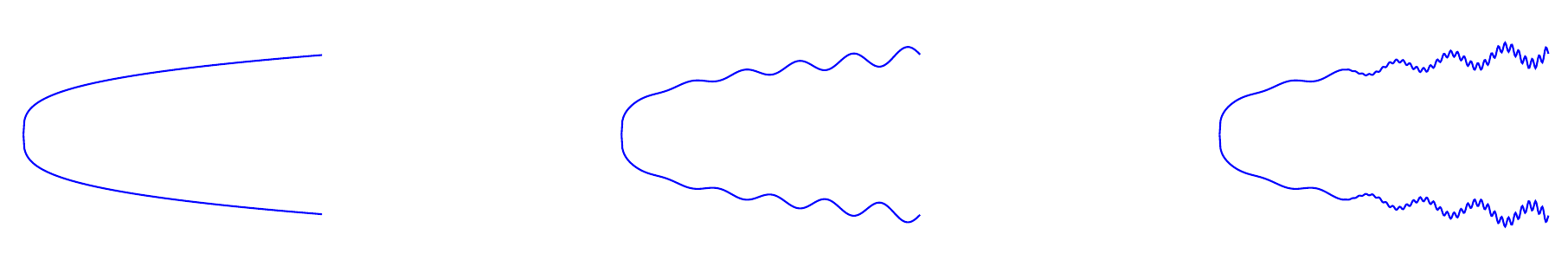}
      \caption{Example of $ \frac{\Om_T}{|\ph (T)|}$ for $T=1$, $T=5$ and $T=20$.}
      \label{Figure 1}
\end{figure}

\fr
The corresponding rescaled map is given by the formula
$$
\ph_T(z):=\frac{1}{|\ph(T)|}\ph(T(z))
$$ 
and one can consider the integral means
$$ 
\b_{\ph_T}(r,t)
:=\frac{\log \int_{1 \leq |y| \leq 2} |\ph_T '(r+iy)|^tdy}{\log 1/r}\,.
$$
Starting from this formula we naturally assign to the tract $\Om$ an integral means spectrum 
$ \b_\infty $
which measures the fractal behavior of the tract at infinity. As in the classical setting, the important function will be
the convex one:
$$
 b_\infty(t):=\b_\infty(t)-t+1\quad , \quad t\in \R\,.
$$
This function always has a smallest zero $\tstar>0$ and, in the good cases, 
$b_\infty$ has a unique zero and is negative in $(\tstar , \infty )$. In this latter case, we will say that the function $f$ has \emph{negative spectrum}.

\medskip

It will become transparent in Proposition \ref{y.4} that there is a strong relation between  the transfer operator, on whose properties the thermodynamic formalisms relies, and the integral means spectrum and we will derive from this that the negative spectrum property implies good behavior of this operator. Negative spectrum turns out to be a very general condition which holds as soon as the tracts have some nice geometry such as the  \emph{H\"older tract} property which essentially means that the domains
$\Om_T$ are uniformly H\"older, see Definition \ref{27} for the precise definition. 
For example, if  the tract itself $\Om$ is a quasidisk then it is a H\"older tract.

\bprop\label{prop neg spectrum}

Let  $f\in \cB$ be an entire function having finitely many tracts.
If the tracts are H\"older then
$f$ has negative spectrum.
\eprop

Throughout the paper we will primarily work with functions of so called disjoint type, which is a particular form of hyperbolicity; we will provide its precise definition in the sequel. We would however like to note that for functions within class $\cS$, by using standard bounded distortion arguments, our results carry over to all hyperbolic functions (see Section \ref{dishyp}) and not merely those of disjoint type. Class $\cD$, for which most of our main results will be formulated and proved, essentially consists in disjoint type functions of class $\cB$ having finitely many tracts, see Definition \ref{9}. 

\bthm\label{thm intro}
Let  $f\in \cD$ be a function having negative spectrum and let $\tstar\in  ]0,2]$  be the smallest zero of $b_\infty$. Then, the following holds:
\begin{itemize}
  \setlength\itemsep{1mm}
\item[-]  For every $t>\tstar$,  the whole thermodynamic formalism, along with its all usual consequences holds: the Perron-Frobenius-Ruelle Theorem, the Spectral Gap property along with its applications: Exponential Mixing, Exponential Decay of Correlations and Central Limit Theorem  (see Section \ref{section MD}).
\item[-]  For every $t<\tstar$, the series defining the transfer operator $\pft$ (see \eqref{29}) is divergent.
\end{itemize}
\ethm

\
Therefore, thermodynamic formalism is crystal clear for functions in class $\cD$ with negative spectrum. The proof 
is based on Theorem \ref{50} which is valid for all functions in class $\cB$ without any further assumptions.

\medskip

There is a very general result of approximating a model function by entire functions due to  Bishop \cite{Bishop-EL-2015, Bishop-S-2016}. His work is
motivated by earlier results of Rempe-Gillen \cite{Rempe-HypDim2}. 
We show in Proposition  \ref{z.10} that the H\"older tract property is preserved when passing from the model
to the approximating entire function. 
In fact, as Lemma~\ref{w.4} demonstrates, the H\"older tract property is  a quasiconformal invariant.
This has a second important application: for entire functions of class $\cS$ the H\"older tract property is 
in fact a property of an analytic family of functions and not only of a single function.
More precisely, if $g\in \cS$  then Eremenko-Lyubich \cite{EL92} naturally associated to $g$ an analytic family of entire functions $\cM_g \subset \cS$. Proposition \ref{w.3} states that every function of $\cM_g$ has H\"older tracts if a function, for example $g$, has. A concrete application of all of this is the following.

\bthm \label{thm 1.4}
Let $g\in \cS$ be any function having finitely many tracts over infinity and assume that they are  H\"older.
 Then every function $f\in \cM_g$ has negative spectrum and the thermodynamic formalism holds for every hyperbolic map from $\cM_g$.
\ethm

We also study a particular family of entire functions called Poincar\'e functions studied previously in \cite{DO08, MP12, ER15} among others. If
$p:\oc\to\oc$ is a polynomial  and if $z_0\in \J_p$ is a repelling fixed point of $p$ then there exists an entire function  $f:\C\to\C$ such that 
\beq\label{linearizer}
f(p'(z_0)z) = p(f(z))
\quad \text{for all $z\in\C$.}
\eeq

\begin{figure}[h]
   \includegraphics[height=7.2cm]{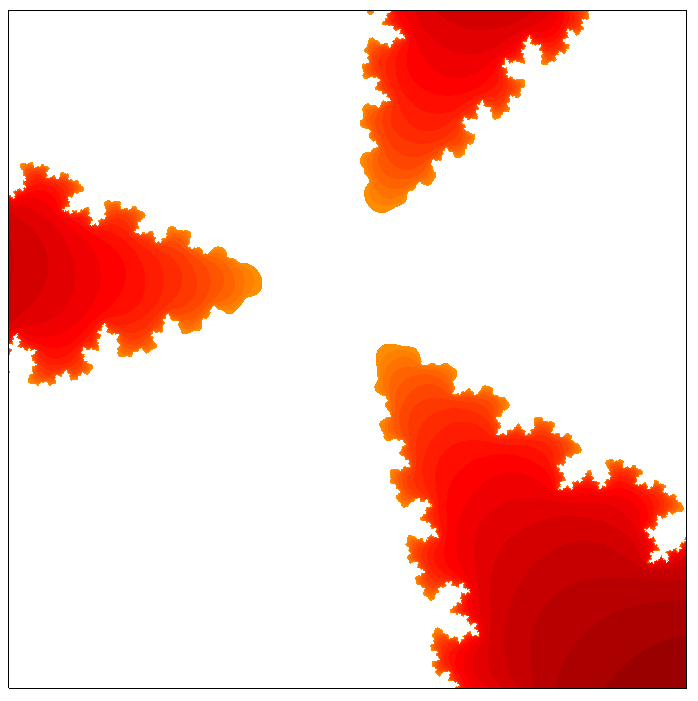}
      \caption{Three fractal and H\"older tracts from a linearizer of Douady's Rabbit.}
      \label{Figure 3}
\end{figure}

\fr
For all entire functions $f$ that obey such a particular linearizing functional equation
such that the involved polynomial $p$ has a connected Julia set
 we show, by a direct calculation in Theorem~\ref{irr pfunctions}, that the transfer operator behaves well. But not all of them have negative spectrum.
Based on the work of Graczyk, Przytycki, Rivera-Letelier and Smirnov \cite{GS98, PRS03}, we show that Poincar\'e functions
 have H\"older tracts if and only if the corresponding linearizing polynomial is topological Collet-Eckmann (TCE).
 In addition, such a linearizer is in  $\cS$ if and only if the polynomial is post-critically finite (thus TCE). 
 Therefore, such functions can be taken 
as generating function of the analytic family in Theorem \ref{thm 1.4}. They are particularly  intriguing since it follows from Zdunik's Theorem~\ref{x.6} in \cite{Zdunik90} that the tracts of Poincar\'e functions are fractals except for the case of polynomials of the form $z\mapsto z^d$, $d\ge 2$, or Tchebychev ones.

 \bcor  \label{thm 1.4'}
 Let $g\in \cS$ be a Poincar\'e function of a polynomial having connected Julia set.
 Then every function $f\in \cM_g$ has negative spectrum and the thermodynamic formalism holds for every hyperbolic map from $\cM_g$.
 \ecor

\smallskip

\sp
This also leads to geometric applications
provided that the topological pressure, as defined in Section \ref{section MD'}, has a zero $h>\tstar$. 
The following result completes the picture on various Bowen's Formulas
(see \cite{MyUrb08} but especially  \cite{BKZ-2012}  which contains a very general version of it).

\bthm[Bowen's Formula]\label{Bowen's}
Let $f\in \cD$ have negative spectrum and be such that the topological pressure $\P(t)$ has a zero $h>\tstar$. Then, the hyperbolic dimension $\HypDim(f)$
of $f$ is equal to the unique zero $h>\tstar$ of the topological pressure.
\ethm

\fr
Here is an other concrete application which is based on the previous formula.

\bthm[Real analyticity of hyperbolic dimension]\label{analyticity}
Let $p:\oc\to\oc$ be a hyperbolic polynomial with connected Julia set, let $z_0\in \J_p$ be a repelling fixed point of $p$
and let $f:\C\to\C$ be an entire function such that \eqref{linearizer} holds. Let the entire functions $f_\ka$ be given by 
$f_\ka(z)=f(\ka z)$, $z\in \C$.
Then, there exists $r>0$ such that  the function
$$
\ka \longmapsto \HypDim(f_\ka ) 
$$
is real analytic in $\D (0,r)\setminus \{0\}$
and $\HypDim(f_\ka ) > \Hdim (\J_p )$.
\ethm

\smallskip
In conclusion, we get a complete, natural and quite elementary approach for the thermodynamic formalism for 
entire functions having negative spectrum which goes far beyond the existing setting since the entire functions in $\cD$ that satisfy the balanced growth condition (see  \cite{MyUrb08, MUmemoirs}) in the tracts  have negative spectrum and, furthermore, they are \emph{elementary} in the sense that the integral means spectrum is as simple as possible: namely $\b_\infty \equiv 0$ and  $\tstar=1$ (see Proposition \ref{17}). 

\ 

\fr{\bf Additional Remark.} 
After having sent out the first version of this paper, Dezotti and Rempe-Gillen informed us that they are actually finishing the preprint \cite{DR} and supplied us with its preliminary version. Concerning thermodynamic formalism, they establish its version for hyperbolic Poincar\'e functions of TCE polynomials. In particular, they show that our Proposition \ref{38} holds for TCE polynomials.

\medskip

\section{The Setting} \label{Section 2}

Let $f:\C\to\C$ be an entire function and let $S(f)$ be the closure of the set of critical values and finite asymptotic values of $f$.
The classification of all types of singularities of an entire function, known as Iversen's classification, is very well explained and presented in \cite{BE08.1}.
We consider functions of the Eremenko--Lyubich class $\cB$ which consists of all entire functions for which the set $S(f)$ is a bounded set. This class contains an important subclass, called Speiser class, which consists of all entire functions for which the set $S(f)$ is finite. 

The dynamical setting is the following. An arbitrary entire function $f:\C\to\C$ is called \emph{hyperbolic} if $f\in \cB$ and if there is a compact set $K$ such that 
$$
f(K)\subset \Int(K)
$$ 
and $f:f^{-1}(\C\setminus K)\to \C\setminus K$ is a covering map. According to Theorem 1.3 in 
\cite{RempeSixsmith16}, an entire function $f$ is hyperbolic if and only if the postsingular set 
$$
\P(f):= \overline { \bigcup_{n\geq 0} f^n (S(f))}
$$
is a compact subset of the Fatou set of $f$. In particular, we have then
\beq\label{hyp dist}
\dist \left( S(f) , \J_f\right) \geq \dist \left( \P(f) , \J_f\right)>0\,.
\eeq
Here and in the sequel, $\J_f$ stands for the Julia set of $f$ defined in the usual way (see for example the survey \cite{BW}).

Concerning the radial Julia set, there are several definitions in the literature (see \cite{MUmemoirs, Rempe09-hyp}).
It is explained in Remark 4.1 of \cite{MUmemoirs} that these definitions lead to different sets whose difference  is dynamically insignificant. In particular they 
have the same Hausdorff dimension. Since we deal only with hyperbolic entire functions, the following definition
fits best to our context:
$$\J_r(f) = \{z\in \J(f) \;: \; \liminf_{n\to\infty} f^n(z) <\infty\}\,.$$
The \emph{hyperbolic dimension} of $f$ is the Hausdorff dimension of this set:
$$\HypDim(f)= {\rm Hdim} ( \J_r (f))\,.$$

Of crucial importance for us is the concept of \emph{disjoint type}. It first implicitly appeared in \cite{B07} and has been explicitly studied in several papers including \cite{RRRS, Rempe09, R-ALC}. In these papers it meant that the compact set $K$ in the definition of a hyperbolic function can be taken to be connected. In this case, the Fatou set of $f$ is connected. We will use its normalized form described below. 

For every $r>0$ let
$
\D_r:=\D(0,r)
$
be the open disk centered at the origin with radius $r$ and 
$
\D_r^*=\C\sms \ov\D_r
$
for the complement of its closure. We denote
$
A(r,R):=\D_R\setminus \ov \D_r
$
 the annulus centered at $0$ with the inner radius $r$ and the outer radius $R$. We further write
$
\D:=\D_1
$
for the unit disk in $\C$ and
$
\D^*:=\D_1^*
$
for the complement of its closure. 

If 
$$
S(f)\subset \D
$$ 
then $f^{-1}(\D^*)$ consists of mutually disjoint unbounded Jordan domains $\Omega$ with real analytic boundaries such that $f:\Om\to \D^*$ is a covering map (see \cite{ EL92}). In terms of the classification of singularities, this means that $f$ has only logarithmic singularities over infinity.
These connected components of $f^{-1}(\D^*)$ are called tracts and 
the restriction of $f$ to any of these tracts $\Om$ has the special form 
\beq\label{11}
\text{$f_{|\Om} =exp\circ{\tau}\;\;$
where $\;\;\ph=\tau^{-1} :\cH:=\{z\in\C:\Re (z) >0\} \to \Om$}
\eeq
is a conformal map. We will always assume that $f$ has only finitely many tracts:
\beq\label{8}
f^{-1}(\D^*)=\bigcup_{j=1}^N \Om_j\,.
\eeq
Notice that this is always the case if the function $f$ has finite order. Indeed, if $f$ has finite order then the Denjoy-Carleman-Ahlfors Theorem (see \cite[p. 313]{Nevbook74}) states that $f$ can have only finitely many direct singularities and so, in particular, only finitely many logarithmic singularities over infinity. 

If $f\in \cB$ is such that 
\beq\label{1}
S(f)\subset \D \quad  \  \text{and} \quad  \  \bigcup_{j=1}^N \ov\Om_j \cap \ov \D  =\emptyset   \quad  \({\rm equivalently:} \  f^{-1}(\ov{\D^*})=\bigcup_{j=1}^N \ov\Om_j \sbt \D^*\), 
\eeq
then we will call $f$ a function of \emph{disjoint type}. 
This is consistent with the disjoint type models in Bishop's paper \cite{Bishop-EL-2015}.
The function $f$ is then indeed of disjoint type in the sense of \cite{B07, RRRS, Rempe09} described above as one can take for $K$ the set $\ov{\D}$. Throughout this paper we will always understand the concept of disjoint type in it more restrictive form of \eqref{1}. 

It is well known (see \cite[p.261]{Rempe09}) that for every $f\in \cB$  the function $\l f$, $\l\in \C^*$, is of disjoint type provided $\l$ is small enough.

\medskip

In our present paper we focus on the following class $\cD$ of entire functions.
 
 \bdfn\label{9}
An entire function $f:\C\to \C$ belongs to class  $\cD$ if the following holds:
  \ben

 \sp \item $f$ has only finitely many tracts, i.e. \eqref{8} holds.

 \sp \item $f$ is of disjoint type in the sense of \eqref{1}; in particular $f$ belongs to class $\cB$.
 
 \sp \item For every tract, the function $\ph$ of \eqref{2_2017_06_28} satisfies the following very general geometric condition: there exists a constant $M\in (0,+\infty)$ such that for every $T\geq 1$ large enough,
\beq\label{(4.2)}
|\ph(\xi )|\leq M |\ph (\xi ')| \quad \text{for all \ $\;\xi , \xi ' \in Q_T\setminus Q_{T/8}$}, 
\eeq
where, we recall, the rectangles $Q_S$ have been defined in \eqref{2}. 
\een
 \edfn

Frequently, only the dynamics of the restriction of $f$ to the union of the tracts will be relevant. We recall 
 from \eqref{11} that such a restriction is given on each component $\Om_j$ by a conformal map $\ph_j: \cH \to \Om_j$ which extends continuously to infinity.

\bdfn\label{tract model}
A  model $(\tau , \Tract )$ is a finite union $\Tract=\bigcup_{j=1}^N \Om_j$ of simply connected unbounded domains $\Om_j$ 
along with
conformal maps $\tau_{|\Om_j}=\tau_j : \Tract_j \to \cH$ such that $\ph_j=\tau_j^{-1}$ extends continuously to infinity:
$$
\text{If } \xi_n\in \cH \text{ with }  \lim_{n\to\infty}|\xi_n|=\infty \quad \text{then}\quad 
 \lim_{n\to\infty}|\ph_j(\xi_n)|=\infty \, .
$$
Associated to $(\tau , \Tract )$ is the model function $f=e^\tau$ and we say that $f\in \cD$
 if $f$ is a disjoint type model
in the sense that  $\;  \ov\Om \cap \ov \D  =\emptyset$.
\edfn

The Julia set of a model $f$ is defined by
$$
\J_f:=\{z\in \D^*\, ; \;\;  f^n(z)\in \D^* \; \; \text{for every}\;\; n\geq 1\}\,.
$$
By Proposition 2.2 in \cite{R-ALC}, this definition coincides with the usual definition of the Julia set in the case of a disjoint type entire function.

Given these definitions, we will write in this paper $f\in \cD$ for either an entire or a model function $f$ having the properties of class $\cD$.
Model functions can be approximated by entire functions of class $\cB$. 
Rempe--Gillen 
\cite[Theorem 1.7]{Rempe-HypDim2} has a very precise result on uniform approximation.
 A weaker notion of approximation
of a model function $f$  by an entire function $g$ is when 
there exists a quasiconformal map $\ph$ of the plane such that 
$f=g\circ \ph$.
Bishop in  \cite{Bishop-EL-2015, Bishop-S-2016} 
has established the existence of such quasiconformal approximations in full generality. In his results $\Om$ can be an arbitrary disjoint union of tracts
and he can approximate by functions in class $\cB$ and even in class $\cS$. We will come back to this in Section \ref{holder tracts} when discussing H\"older tracts.

\smallskip
If $f$ is of disjoint type, either entire or model, then
$$
\jul_f  \subset  \bigcup_{j=1}^N\Om_j
$$
and the Julia set $\jul_f $ is entirely determined by the dynamics of $f$ in the tracts. So, for disjoint type functions we can work indifferently either with a model or a global entire function. It follows from \eqref{8} and \eqref{1}
that for such functions
\beq\label{1_2017_06_30}
f^{-1}\Big(\bigcup_{j=1}^N \ov\Om_j\Big)\sbt \bigcup_{j=1}^N\Om_j,
\eeq
and that there exists  $\ga \in (0,1)$ such that 
\beq\label{33}
\J_f \subset \Om= \bigcup_{j=1}^N\Om_j \sbt  \D_{e^\ga}^*.
\eeq
As said, throughout the whole paper we restrict our attention to the functions in class $\cD$, so, in particular, to those of disjoint type. One can extend all our considerations and results to the case of hyperbolic entire functions belonging to Speiser class $\cS$, i.e. replacing \eqref{1} by mere hyperbolicity in Definition~\ref{9} and assuming class $\cS$. This is, quite easily, done in Section~\ref{dishyp} by using  Koebe's Distortion Theorem only (in addition to all what we did for disjoint type functions).

\sp Here and in the sequel we use the classical notation such as 
$$
A \asymp B.
$$ 
As usually, it means that the ratio $A/B$ is bounded below and above by strictly positive and finite constants
that do not depend on the parameters involved. The corresponding inequalities up to a multiplicative constant are 
denoted by 
$$
A\preceq B \  \  \  {\rm and } \  \   \   A\succeq B.
$$
With this notation we have the following. We recall that the rectangle
$Q_2$ has been defined in \eqref{2}.

\bthm[Bounded distortion]\label{bdd distortion}
If $\ph:Q_2\to\C$ is a univalent holomorphic map then, for every $ 0<r<1 $ and every  $-2\leq y\leq 2$, we have that
\beq\label{bdd 1}
|\ph' (1)|(1-r)\preceq |\ph'((1\pm r)+iy)|\preceq |\ph' (1)|\frac{1}{(1-r)^3} \; . 
\eeq
If $\ph$ is a univalent holomorphic map defined on the entire half-plane $\cH$ then, for every $x>1$,
\beq\label{bdd 2}
|\ph' (1)|\frac{1}{x^3} \preceq |\ph'(x)|\preceq |\ph' (1)| x \; . 
\eeq
Here, the multiplicative constants involved are absolute.
\ethm

\bpf
This is simply a fairly straightforward application of Koebe's Distortion Theorem. Let $g:Q_2\to \D$ be conformal, i.e univalent and holomorphic surjection. It has a holomorphic extension
to a neighborhood of $Q_1$ in $\C$ and thus $g_{|Q_1}: Q_1\to g(Q_1 )$ is a bi--Lipschitz map. It suffices thus to apply 
Theorem 1.3 in \cite{PommerenkeBook} to $\ph\circ g^{-1}$ in order to deduce \eqref{bdd 1}.
The inequalities in \eqref{bdd 2} also follow since for every univalent map $\ph$ on $\cH$ the map
$z\mapsto \ph (xz)$ is a univalent map
on $Q_2$ and one can apply \eqref{bdd 1}.
\epf

Let $\ph : \cH \to \Om$ be a conformal homeomorphism. Then  \eqref{bdd 2} implies for every $T\ge 1$ that
\beq\label{bdd 3} 
|\ph (T)-\ph (1)| \leq \int _1^T |\ph' (x)| dx\lek |\ph '(1)| \, T^2\,.
\eeq

\smallskip


\section{Fractal behavior of $\partial \Om$ at infinity}\label{section 3}
We first analyze what happens for one single tract. So, we consider a model $(\tau ,\Om)$ with $\Omega$ a simply connected domain.
In Figure \ref{Figure 1} we illustrated the possible fractal behavior of a tract $\Om$ near infinity
by considering rescalings of the exhaustion domains $\Om_T$ of $\Om$.
Associated to these rescaled domains are the rescaled conformal maps $
\ph_T:\cH\longrightarrow \C $
given by the formula
\beq\label{3''}
\ph_T(z):=\frac{1}{|\ph(T)|} \, \ph(Tz) \quad , \quad z\in \cH.
\eeq
We will frequently treat the maps $\ph_T$ as restricted to the set $Q_2$ and will use the same symbol $\ph_T$ for this restriction. In symbols, we will consider the maps
\beq\label{3}
\ph_T: Q_2\longrightarrow  \frac{1}{|\ph(T)|}\Om_{2T}\quad , \quad T\geq \ga 
\eeq
where, as always, $\ga$ comes from  \eqref{33}.
In particular
\beq\label{3B}
|\ph_T(1)|=1.
\eeq
We denote by $\cF_\Om$ the family of all the functions $\ph_T$, $T\geq \ga$. Since asymptotic properties of this family 
will be crucial, we now make some elementary observations. Let us recall here that we always work under the standard assumption \eqref{(4.2)}.

\blem\label{14} Suppose \eqref{(4.2)} holds. Then, $\cF_\Om$ is a normal family, in the sense of Montel, on $Q_1\setminus Q_{1/8}$ and furthermore
\beq\label{z.15}
\text{ $\frac{1}{T^4}\preceq |\ph_T'(1)|\preceq 1 \;\; $ , $\;\; T\geq \ga$,}
\eeq
 \elem

\bpf
It follows from \eqref{(4.2)} and \eqref{3B} that for every $T\geq \ga$ it holds
\beq\label{1_2017_06_28}
\ph_T(Q_1\setminus Q_{1/8} ) \subset \D(0,M).
\eeq
Normality of $\cF_\Om$ follows thus directly from Montel's Theorem. The left hand side of  \eqref{z.15} is a straightforward consequence of the left hand side of item \eqref{bdd 2} of the distortion Theorem~ \ref{bdd distortion} along with \eqref{bdd 3}, both applied to the map $\ph:\cH\to\C$. Indeed, using them, we get
$$
|\ph_T'(1)|
=\frac{|\ph'(T)|}{|\ph(T)|}T
\gek |\ph'(1)|\frac{T}{T^3|\ph(T)|}
=|\ph'(1)|\frac{1}{T^2|\ph(T)|}
\gek \frac{1}{T^4},
$$
with \eqref{bdd 3} invoked for the last inequality sign. Since $Q_1\setminus Q_{1/8}\spt \D(1,1/2)$, it follows from \eqref{1_2017_06_28} that $\ph_T(\D(1,1/2))\sbt \D(0,M)$. But by Koebe's $\frac14$--Distortion Theorem, $\ph_T(\D(1,1/2))\spt \D(\ph_T(1),|\ph_T'(1)|/8)$. Therefore, $|\ph_T'(1)|\le 8M$, formula \eqref{z.15} is proved. 
\epf

Information of the boundary of the image domain can be obtained by considering integral means spectrum
(see  \cite{Makarov98} and \cite{PommerenkeBook} for the classical case which concerns conformal mappings defined on the unit disk). 
In order to do so, let $h:Q_2\to U$ be a conformal map onto a bounded domain $U$ and define
\beq\label{4}
\b_h(r,t):= \frac{\log \int_I |h'(r+iy)|^tdy}{\log 1/r} \; \text{ ,} \  \ r\in (0, 1) \ \text{ and } \  t\in \R \,.
\eeq
The integral is taken over $I=[-2,-1]\cup[1,2]$ since this corresponds to the part of the boundary of $U$ that is important for our purposes.
\begin{figure}[h]
   \includegraphics[height=2cm]{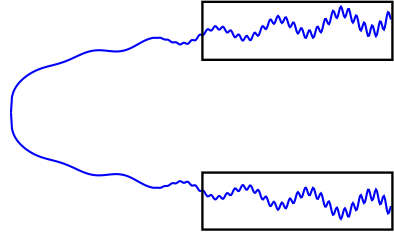}
     \caption{The part of the boundary in the boxes can resemble more and more a fractal as $T\to\infty$.}
      \label{Figure 2}
\end{figure}

\fr A well known application of the distortion Theorem~\ref{bdd distortion} shows that there exists $K_1>0$ such that, for all $t\geq0$,
\beq\label{z.16}
 -t+\frac{t\log(|h'(1)|)-K_1}{\log 1/r} 
\leq \b_h(r,t)
\leq 3t+\frac{t\log\big( |h'(1)|)+K_1}{\log 1/r}\;
\eeq
and a corresponding inequality holds for all  $t<0$.
Replacing now $h$ by the conformal maps of the family $\cF_\Om$, one first has the following observation.

\blem\label{R-independant}
Let $t\in \R$, $R\geq \ga$ and set $T= R/r$, $r>0$. Then
 $$ \limsup_{r\to 0} \b_{\ph_{R/r}}(r,t) = \limsup_{T\to +\infty} \b_{\ph_T}(1/T,t)$$ is finite and does not dependent on $R\geq \ga$.
\elem

\bpf It suffices to treat the case $t>0$ since $t<0$ can be treated the same way and for $t=0$  there is nothing to show. 
So, let $t>0$ and let $R\geq \ga$. 
Finiteness of $ \limsup_{r\to 0} \b_{\ph_{R/r}}(r,t)$ directly results from  \eqref{z.15} and \eqref{z.16}.
In order to study the dependence on $R$ of this expression, observe that
$$
\b_{\ph _{R/r}}(r,t)= \frac{\log \int_I |\ph_{R/r}'(r+iy)|^tdy}{\log 1/r} = 
\frac{\log \int_I \left( \frac{R/r}{|\ph (R/r)|}\right)^t|\ph'(R+i\frac{R}{r}y)|^tdy}{\log 1/r} \,.
$$
Since $0<\ga <1$, it suffices to compare this with the corresponding expression with $R'=1$.
If $r'=r/R$ then
$$
\b_{\ph _{1/r'}}(r',t)= \frac{\log \int_I |\ph_{1/r'}'(r'+iy)|^tdy}{\log 1/r'} = 
\frac{\log \int_I \left( \frac{R/r}{|\ph (R/r)|}\right)^t|\ph'(1+i\frac{R}{r}y)|^tdy}{\log 1/r'} \,.
$$
By Theorem \ref{bdd distortion} there exists $K\geq 1$ such that
$$
\frac1{KR^3}\left|\ph'\left(1+i\tilde y\right)\right|\leq \left|\ph'\left(R+i\tilde y\right)\right|\leq\left|\ph'\left(1+i\tilde y\right)\right| KR
$$
for every $\tilde y \in \R$ and every $R\geq 1$. For $\ga \leq R <1$ a corresponding estimation holds, we omit this detail and consider  $R\geq 1$. Then
$$
\frac{-t\log ( KR^3)}{\log 1/r }+ \b_{\ph _{1/r'}}(r',t)\frac{\log 1/r'}{\log 1/r}\leq
\b_{\ph _{R/r}}(r,t)\leq
\frac{t\log ( KR)}{\log 1/r }+ \b_{\ph _{1/r'}}(r',t)\frac{\log 1/r'}{\log 1/r}
$$
from which follows that 
$$
\limsup_{r\to 0} \b_{\ph_{R/r}}(r,t)  = \limsup_{r'=\frac{r}{R}\to 0} \b_{\ph_{1/r'}}(r',t)\,. 
$$
\epf

\fr
We are now ready to introduce the  function
\beq\label{5a}
\b_\infty (t) = \limsup_{r\to 0} \b_{\ph_{1/r}}(r,t)= \limsup_{T\to +\infty} \b_{\ph_T}(1/T,t)\quad , \quad t\in \R\,.
\eeq
Lemma \ref{R-independant} justifies that $t\mapsto \b_\infty (t)$ is a well defined finite function on $\R$ and that
\beq\label{5b}
\b_\infty (t) = \sup_{R\geq \ga }\;  \limsup_{r\to 0} \b_{\ph_{R/r}}(r,t) \quad , \quad t\in \R \,.
\eeq
Up to now we considered a single tract.
In the general case we deal with a function $f\in \cD$ and so $f^{-1}(\D^*)$ is a disjoint union of finitely many
tracts $\Om_j$, $j=1,...,N$. Denoting  $\b_{\infty ,j}(t)$ the function of \eqref{5a} defined in the tract  $\Om_j$, we can associate to $f$ the function
$$
\b_\infty=\b_{\infty ,f} :=\max_{j=1,...,N}  \b _{\infty ,j}\,.
$$

\fr
Now we continue dealing with one fixed tract and we skip the index $j$.
\medskip

\bprop\label{6} 
The function $ \b_\infty :\R\to \R$ is convex with
$$
\b_\infty (0)=0 \  \  \  \text{and } \  \  \b_\infty (2)\leq 1.
$$
\eprop

\bpf
All involved $\b$--functions are convex by a classical application of H\"older's inequality (see for example p.176 in \cite{PommerenkeBook}). It is trivially obvious that $\b_\infty(0)=0$ while $\b_\infty (2)\le 1$ results from the well known area estimate. Indeed, let $0<r<1$. For every integer $0\leq k < 1/r$ set
$
y_k^+:= 1 +kr$,  
$$
U_k^+:=\big\{z\in\C: r<\Re z <2r \  \ {\rm and } \  \  y_k^+ <\Im z <y_{k+1}^+\big\}
\quad \text{and} \quad 
U_k^-:=\big\{ \ov z \; , \, z\in U_k^+\big\} \; .
$$
Then,
$$
\int_I |\ph_T'(r+iy)|^2 dy \asymp \sum_{\substack{0\leq k <[1/r]  \\ \ep \in \{ + , -\}}} \frac1r area(\ph_T(U_k^\ep))
\leq \frac1r area (\ph_T (Q_1\setminus Q_{1/8})) \preceq \frac1 r
$$
since $\diam\(\ph_T (Q_1\setminus Q_{1/8})\) \lek 1$. Applying logarithms, dividing by $\log(1/r)$ and letting $r\to 0$ gives that $\b_\infty (2)\leq 1$. 
\epf

\fr
A function related to $\b_\infty $, that will be crucial in the sequel, is the following:
\beq\label{16}
 b_\infty (t):=\b_\infty (t)-t+1 \quad , \quad t\in \R\,.
\eeq
As an immediate consequence of Proposition~\ref{6} we get the following.

\bprop\label{p1_2017_07_10}
The function $b_\infty:\R\to\R$ is also convex, thus continuous, with 
$$
b_\infty (0)=1 \  \ {\rm and }  \  \  b_\infty (2)\leq 0.
$$
\eprop

\medskip

\fr Consequently, the function $b_\infty$ has at least one zero in $]0,2]$ and we can introduce a number 
$\tstar\in (0,2]$ by
\beq\label{7}
\tstar:= \inf \{ t>0 \; : \; b_\infty (t)=0\}
= \inf \{ t>0 \; : \; b_\infty (t)\le 0\}\,.
\eeq
Again, in the case of a function $f$ with finitely many tracts $\Om_j$ we thus have finitely many numbers 
$\tstarj$ and then we set 
$$
\tstar:= \max_j \tstarj.
$$

We will
consider below various situations and examples illustrating the behavior of $b_\infty$ and of $\tstar$. 
Notice also that the paper \cite{My-HypDim} also is based on  $\b_\infty$ 
along with the zero (called also $\Theta$) of $b_\infty$.

\

In order to perform full thermodynamic formalism we need the following crucial property.

\bdfn\label{18} 
A function $f\in \cD$ has negative spectrum if
$$ b_\infty (t)<0 \; \text{ for all} \ \ t> \tstar\,.$$
\edfn

\fr
As we will see in Section \ref{sec reg funct}, this property does hold if the tracts have some nice geometry.

\medskip

\section{Transfer operator}
In the sequel  $f$ will be either an entire function in $\cD$ or a model map in $\cD$ and
we will work with the Riemannian metric 
\beq\label{3_2017_06_28}
|dz|/|z|.
\eeq
This metric is conformally equivalent with the standard Euclidean one, it has singularity at $0$ but is tailor crafted for our analysis of Perron--Frobenius operators. 
The derivative of a holomorphic function $h$ calculated with respect to the metric of \eqref{3_2017_06_28} at a point $z$ in the domain of $h$ is denoted by $|h'(z)|_1$ and is given by the formula
\beq\label{3_2017_06_28bis}
|h'(z)|_1=|h'(z)|\frac{|z|}{|h(z)|}.
\eeq

So, given a real number $t\geq 0$, we define the transfer operator $\pft$ by the usual formula:
\beq\label{29}
\pft g (w):= \sum_{f(z)=w} |f'(z)|^{-t}_1 g(z) \quad \text{for every}\quad w\in \ov\Om\,
\eeq
where $g$ is any function in $\cC_b (\ov\Om)$, the vector space of all continuous bounded functions defined on $\Om$. The norm on this space, making it a Banach space, will be the usual sup-norm $\|\cdot\|_\infty$.
 Note that if $w\in\ov\Om$, then $f^{-1}(w)\sbt\Om$, whence $f'(z)$ is well defined for all $z\in f^{-1}(w)$ and, in consequence, all terms of the above series are also well defined.

\bthm \label{50}
Let $f$ be a model or an entire function of class $\cB$ such that $S(f)\subset \D$.
Assume that there exists $s>0$ and $w_0\in \D^*$ such that 
$$
\pft \1(w_0)<\infty \quad \text{for every} \quad  t>s\,.
$$
Let $\tilde\ga \in (0,1)$. Then
$$
 \sup_{|w| > e^{\tilde\ga} }   \pft \1 (w) <\infty \quad \text{for every} \quad  t>s \, .
$$
In addition,  for all $t>s$ and $p>1$ such that $\frac1p< \frac{t}{s}-1$, there exists a constant $C_{p,t}$ such that
\beq\label{55}
\pft\1 (w) \leq  \frac{C_{p,t}}{(\log|w|)^{1/p}} \quad \text{for all} \quad w\in \D_{e^{\tilde\ga}}^*\,.
\eeq
\ethm

\fr In this result, no dynamical hypothesis nor finiteness of the number of tracts is assumed. If we restrict
to functions where $\Omega$ is backward invariant then it tells us that the transfer operators $\pft$ are bounded.

\bcor \label{50'}
Let $f\in \cD$. If there exists $s>0$ and  $w_0\in \Om$ such that $\pft \1(w_0 )<\infty$ for all $t>s$, then all $\pft$, $t>s$, are bounded operators of $\cC_b (\ov\Om)$, satisfying in addition \eqref{55}.
\ecor

We will explain in Section \ref{section MD} that the conclusion of this result combined with our previous work \cite{MUmemoirs}
leads to full thermodynamic formalism along with all its usual consequences.

\bpf Although $f$ may have infinitely many tracts, it suffices to consider the case of a single tract $\Omega$ since the estimates we obtain generalize directly. 
If $w\in\ov\D_{e^{\tilde\ga}}^*$ then
\beq\label{3_2017_06_30}
\pft \1 (w) = \sum_{f(z)=w} |f'(z)|^{-t}_1 \,.
\eeq
The function $f$ restricted to $\Om$ is of the form $f=e^\tau$. Thus
$$
|f'(z)|_1= \frac{|f'(z)|}{|f(z)|}|z| = |\tau '(z)||z|\,.
$$
Since $\ph =\tau ^{-1}$, we have that
$$
|f'(z)|_1 = \left| \frac{\ph (\xi)}{\ph'(\xi)}\right|,
$$ 
where $\xi = \tau (z)$. In the series of \eqref{3_2017_06_30} $z$ runs through the preimages of $w$ under $f$, thus $\xi$ runs through the set $\exp^{-1}(w)$. Let 
$$
R:=\log |w|=\Re (\xi)
$$
for every $\xi\in\exp^{-1}(w)$. We have $R\geq\tilde\ga >0$. We have 
\beq\label{51}
\pft \1 (w) = \sum_{\xi\in\exp^{-1}(w)} \left| \frac{\ph' (\xi)}{\ph(\xi)}\right|^t = \sum_{\xi\in\exp^{-1}(w)} \left| (\log \ph)'(\xi)\right|^t 
\eeq
with an arbitrary choice of a holomorphic branch of the logarithm of $\ph$. Koebe's Distortion Theorem applied to the conformal map
$\log \ph : \cH \to \log \Om$ gives
\beq\label{53}
\pft \1 (w) \asymp \int _\R  \left| (\log \ph)'(R+iy)\right|^t dy
\eeq
On the other hand, $\log \ph$ is an inverse branch of the logarithmic coordinates of the function $f$ 
as defined in Section 2 of \cite{EL92}. Hence, Lemma 1 of \cite{EL92} applies and yields
$$
 \left| (\log \ph)'(\xi)\right| \leq \frac{4\pi}{\Re \, \xi } \quad , \quad \xi \in \cH\,.
$$
In particular, the holomorphic function $u:\cH\to \C$, defined by 
$$
u(z):=(\log \ph)'(\frac{\tilde\ga}{2} +z),
$$
is bounded
and the function $z\mapsto |u(z)|^t$ is subharmonic, continuous on $\ov \cH$ and bounded. We can therefore compare it 
with its harmonic majorant as it is done in  \cite[Corollary 10.15]{Mash09}:
\beq\label{52}
| u(z)|^t \leq \frac1\pi \int _\R \frac{x}{(y-s)^2+x^2}|u(is)|^tds \quad , \quad z=x+iy \in \cH\, . 
\eeq
Integrating this inequality and using Fubini's Theorem gives
\beq\label{54}
\int_\R | u(x+iy)|^t \, dy \leq \int _\R \frac1\pi \int_\R \frac{x}{(y-s)^2+x^2}dy \,  |u(is)|^tds =  \int _\R |u(is)|^tds <+\infty
\eeq
where the last inequality holds since, by Koebe's Distortion Theorem and  \eqref{53}, the integral on the right hand side is comparable to 
$\pft\1 (w_0)$, where  $w_0\in \Om$ is the point of the assumptions in Theorem \ref{50}.
Therefore,  \eqref{53} along with \eqref{54} imply that
$$
\pft \1 (w) \preceq \pft\1(w_0) \quad \text{for every}\quad w\in \D_{e^{\tilde\ga}}^*.
$$
We have thus proved that $\pft$ is uniformly bounded on $\ov\Om$.

It remains to show the additional property \eqref{55}. Write $t=\tau+\delta$ where $\tau > s$ and  $0< \d < t-s$.
With $x=\log|w| -\tilde \g /2$, formulas \eqref{53} and \eqref{54} imply that
$$
\pft\1 (w) \preceq \frac1\pi \int _\R  \int_\R \frac{x}{(y-s)^2+x^2} \,  |u(is)|^\d ds \,  |u(x+iy)|^\tau dy\, .
$$
For every $p,q>1$ such that $\frac1p+\frac1q=1$ and $p\d > \tstar$, H\"older's inequality yields
$$
 \int_\R \frac{x}{(y-s)^2+x^2} \,  |u(is)|^\d ds\leq \left(  \int_\R  \left( \frac{x}{(y-s)^2+x^2} \right)^qds  \right)^\frac1q
 \left( \int _\R |u(is)|^{\d p} ds \right)^\frac1p =O (x^{-\frac1p})\,.
$$
Combining the last two displayed formulas, we get
$$
\pft\1 (w) \leq \frac{C_1}{x^{1/p}} \int _\R   |u(x+iy)|^\tau dy
\leq \frac{C_2}{x^{1/p}} \quad \text{for every} \quad w\in \D_{e^{\tilde \ga}}^*
$$
for appropriate constants $C_1,C_2$ depending on $\d$ and $p$. The proof is now complete since $x= \log |w| -\ga /2\asymp \log |w|$
for every $w\in \D_{e^{\tilde \ga}}^*$.
\epf

\smallskip

Given Theorem \ref{50}, the essential question is to decide whether, for a given function $f\in \cD$ and parameter $t$, the transfer operator evaluated at $\1$ is finite or not at some point. This is where our new geometric tools come into play. Aiming to prove Theorem \ref{thm intro} we first reformulate $\pft$ in terms of the $\b$--functions.

%

\bprop\label{y.4}
If $f\in \cD$ and  $t\geq 0$, then
$$
\pft \1 (w) \asymp (\log |w|)^{1-t} \left\{\int_{-1}^1\left|\ph_{\log|w|} '(1+iy)\right|^t  dy+ 
\sum_{n\geq 1} 2^{n \big(1-t+\b_{\ph_{2^n\log|w|}}( 2^{-n},\, t )\big)}\right\}
$$
for every $w\in\Om$ with the above series being possibly divergent. 
\eprop

\fr The issue of convergence of the above mentioned series will be the next step.

\bpf [Proof of Proposition~\ref{y.4}]
From the proof of Theorem \ref{50} we already have the reformulation of the transfer operator 
that we need in \eqref{51}. Hence
$$
\pft \1 (w) = \sum_{\xi\in \exp^{-1}(w)} \left|\frac{\ph '(\xi)}{\ph (\xi )}\right|^t
=\sum_{\xi\in \exp^{-1}(w)|\atop|\Im \xi| <R}\left|\frac{\ph '(\xi)}{\ph (\xi )}\right|^t + \sum_{n\geq 0} \; \sum_{\xi\in \exp^{-1}(w)\atop 2^nR\leq |\Im \xi| <2^{n+1}R}\left|\frac{\ph '(\xi)}{\ph (\xi )}\right|^t\,.
$$
Applying to each of these sums \eqref{(4.2)} respectively with $T=R$ and $T=2^{n+1}$, $n\ge0$, we get
\beq\label{11_2017_08_05}
\pft \1 (w) \asymp 
\frac{1}{|\ph (R)|^t} \sum_{\xi\in \exp^{-1}(w)|\atop|\Im \xi| <R}\left|\ph '(\xi)\right|^t + 
\sum_{n\geq 0} \frac{1}{|\ph (2^{n+1}R)|^t}\sum_{\xi\in \exp^{-1}(w)|\atop2^nR\leq |\Im \xi| <2^{n+1}R}\left|\ph '(\xi)\right|^t\,.
\eeq
Since two consecutive elements of $\exp^{-1}(w)$ are at distance $2\pi$ and since $\Re \xi = R \geq \ga >0$, Koebe's  Distortion Theorem yields
$$
\pft \1 (w) \asymp \frac{1}{|\ph (R)|^t} \int_{-R}^R\left|\ph '(R+iy)\right|^t  dy+ 
\sum_{n\geq 0} \frac{1}{|\ph (2^{n+1}R)|^t}\int_{ I_{n,R}} \left|\ph '(R+iy)\right|^t dy
$$
where 
$$
I_{n,R}= [-2^{n+1}R, -2^nR]\cup[2^nR, 2^{n+1}R].
$$
Remember that $I=[-2,-1]\cup[1,2]$ and that we have introduced the rescaled functions $\ph_R$ in \eqref{3}. A change of variables gives now
$$
\pft \1 (w) \asymp \; R^{1-t} \int_{-1}^1\left|\ph_R '(1+iy)\right|^t dy+ 
\sum_{n\geq 0} \big(2^{n+1}R\big)^{1-t} \int _I  \left|\ph_{2^{n+1}R} '\left(\frac{1}{2^{n+1}}+iy\right)\right|^t dy\,.
$$
With invoking the definition \eqref{4} this completes the proof of Proposition~\ref{y.4}.
\epf

Passing to functions having negative spectrum, we can now fully describe the behavior of their transfer operators.
As it will be explained in Section \ref{section MD}, this then allows us to prove Theorem \ref{thm intro}, and its usual consequences, following \cite{MUmemoirs}. We recall that for functions with negative spectrum $\tstar$ is the unique zero of $b_\infty$.

\bthm\label{56}
If $f\in \cD$ is a function with negative spectrum, then:
\begin{itemize}
\setlength\itemsep{1mm}
\item[-]  For every $t>\tstar$,  $\| \pft \1 \|_\infty <+\infty$ and \eqref{55} holds.
\item[-]  For every $t<\tstar$, the series defining $\pft \1$ is divergent at every point.
\end{itemize}
\ethm

\bpf 
Let $\ga $ be the constant from \eqref{33}, let $w_0\in \D^*_{e^\ga}$ be any point  and set $R=\log |w_0| >\ga$.
Since $f$ has negative spectrum,  $b_\infty (t)=-2a_t<0$ for every $t>\tstar$. 
It thus follows right from the definition of $\b_\infty$ in \eqref{5a} that
there exist $n_{R,t}>0$ such that
$$
\b_{\ph_{2^nR}}(2^{-n}, \,t)-t+1
\leq -a_t<0 \quad \text{for all} \quad n\geq n_{R,t}\,.
$$
Applying Proposition \ref{y.4}, we get that 
$$
\pft\1 (w_0) <\infty
$$ 
for all $t>\tstar$. We therefore have checked the hypotheses of Theorem \ref{50}. It implies that Theorem \ref{56} holds for all $t>\tstar$.

Let now $t<\tstar$, and $w\in \Om$ be any point. Set again 
$$
R:=\log |w|.
$$
Then $b_\infty (t) =4a_t>0$ and thus, by Lemma \ref{R-independant},  
$$
\limsup_{r\to 0} \b_{\ph _{R/r}}(r,t) -t+1=4a_t>0\; .
$$
Fix a sequence $r_j\downto 0$ such that $\b_{\ph _{R/r_j}}(r_j,t)-t+1\geq 2a_t$ for all $j\ge 1$. Then associate to every  $j\geq 0$ an integer $n_j$ such that
$2^{-n_j-1}<r_j\leq 2^{-n_j}$. Writing down the definition of $\b_{\ph_T}(r,t)$ and employing \eqref{(4.2)} along with bounded distortion, one gets
$$
\lim_{j\to \infty }\;  \frac{\b_{\ph _{R/r_j}}(r_j,t) }{ \b_{\ph _{R2^{n_j}}}(2^{-n_j},t)}=1\,.
$$
Thus, $\b_{\ph _{R2^{n_j}}}(2^{-n_j},t)-t+1\geq a_t$ for all sufficiently large $j$. This implies that the coefficients in the series in
Proposition \ref{y.4} do not converge to zero, whence $\pft\1 (w)=\infty$.
\epf

\smallskip

\section{Functions with negative spectrum and H\"older tracts.} \label{sec reg funct}

Theorem~\ref{thm intro} decisively shows that the transfer operators $\pft$ of negative spectrum entire functions behave 
sufficiently well so that a fairly complete account of the corresponding thermodynamic formalism can be derived. In the current section we want to get some idea of which functions in class $\cD$ may have  negative spectrum. We will start with considering some classical examples such as exponential functions and we will see that the class of balanced functions
in \cite{MyUrb08, MUmemoirs}
behaves like these classical examples (Proposition \ref{17}); it has the simplest possible $b_\infty$ spectrum, namely $b_\infty (t)=1-t$. Functions with such spectrum will be called \emph{elementary}.

Then we will show that a function has negative spectrum as soon as its tracts have some nice geometry. For us it will be H\"older domains. Particular examples of such tracts are quasidisks or tracts having the John or H\"older property used in \cite{My-HypDim}. 
We finally show that functions of infinite order can also have negative spectrum, and thus the thermodynamic formalism applies to them too.

 \subsection{Classical functions and balanced growth}
 
 The most classical transcendental family is certainly $\l e^z$ or, more generally, $\l e^{z^d}$, $\l\in\C\sms\{0\}$, $d\geq 1$.
 By a straightforward calculation (see also Proposition \ref{17}), all these functions 
 have a trivial integral means spectrum $\b_\infty\equiv 0$ and thus 
 \beq \label{25} b_\infty (t)=1-t\quad , \quad t\in \R\,.\eeq
 In particular, they have negative spectrum with $\tstar=1$ and the tracts of  $\l e^{z^d}$ are not fractal at all. 
 This is also clear when we consider the rescalings. For any tract of such a function, the part of its boundary depicted in the boxes in Figure~\ref{Figure 2} converges to a straight line segment as $T\to\infty$.
 
\medskip
 
The thermodynamic formalism has been for the first time developed for some transcendental  meromorphic functions by Krzysztof Bara\'nski in \cite{Baranski95}. He did it for for the tangent family. Then this theory has been established for several other families of meromorphic functions. One should mention here quite a large and general class of meromorphic functions considered in \cite{KU02} where, however, as in \cite{Baranski95}, there were no singular values of $f^{-1}$ in the Julia which was a compact subset of $\oc$. We also mention \cite{UZ03, UZ04}, where for the first time the thermodynamic formalism was built for a transcendental meromorphic function having such singularity in the Julia set. It was in fact $\infty$ being the asymptotic value of hyperbolic exponential functions. The most general, actually the only general, framework comprising all the classes mentioned above and much more, for which a full fledged thermodynamic formalism has been developed, is up to now the one of \cite{MyUrb08, MUmemoirs}. Indeed, these two works cover many classes of entire and meromorphic functions that include such classical functions as exponential family, the ones of the sine and cosine-root family, elliptic functions and all the functions having polynomial Schwarzian derivative.
It is based on a condition for the derivative which, for entire functions takes on the following form.

 \bdfn\label{23}
An entire function $f:\C\to\C$ is said to be of balanced growth if it has finite order, denoted in the sequel by $\rho=\rho(f)$, and if 
 \beq\label{24}
 |f'(z)| \asymp |f(z)|\, |z|^{\rho -1} \quad \text{,} \quad z\in \J_f\,.
 \eeq
 \edfn

 The examples in \cite{MyUrb08, MUmemoirs} that satisfy this condition have non-fractal tracts precisely as the classical exponential functions $\l e^{z^d}$. This is a general fact for balanced functions. They are elementary in the sense that their integral means spectrum $\b_\infty$ is the most trivial possible. In the next result we 
 assume balanced growth on the whole tracts, a condition that is satisfied by
all the examples in \cite{MyUrb08, MUmemoirs}.
 
 \bprop\label{17} If $f\in \cD$ satisfies  the balanced condition \eqref{24} in $\Om$, then $f$ is elementary in the sense that $\b_\infty \equiv 0$ and thus
  $$b_\infty (t)= \b_\infty (t)-t+1= 1-t \quad , \quad t\in \R\,.$$
 In particular, $f$ has negative spectrum with $\tstar=1$.
  \eprop
  
 \bpf It suffices to consider an entire function or a model function $f\in \cD$ with $\Om$ being one single tract. Then
$$
f_{|\Om}=e^{\tau},
$$ 
where 
$$
\tau=\ph^{-1} : \Om \to \cH
$$
is a conformal map. Shrinking $\Om$ if necessary, we may assume that $\tau $ is a continuous map defined on $\ov \Om$.
 We also may assume that $0\not\in\Om$ and that a holomorphic branch of $z\mapsto z^\rho$ can be well defined on $\Om$, where $\rho$ comes from
 \eqref{24}. This allows us to introduce a map 
$$
 h:=\ph^\rho : \cH \to \Om^\rho :=\{z^\rho \; : \; z\in \Om \}\,.
$$
By assumption, $f$ satisfies  \eqref{24}
and $f'=f\, \tau'$.
Therefore, 
$$
|\ph'(\xi)|=\frac{1}{|\tau'(\ph(\xi))|}\asymp |\ph(\xi)|^{1-\rho},
$$
which implies
\beq\label{x.11}
|h'(\xi )| = \rho | \ph (\xi )|^{\rho-1} |\ph' (\xi ) | \asymp 1 \quad , \quad \xi \in \cH\,.
\eeq
As immediate consequence we get that $|h(T)-h(0)|
\preceq T$ which implies 
\beq\label{x.35}|\ph (T)|\preceq T^{1/\rho} \quad \text{for}\quad T\geq T_0\eeq
where $T_0\geq 1$ is such that $|h(T)|\geq 2|h(0)|$ for $T\geq T_0$; it is finite due to \eqref{x.11} and Koebe's $\frac14$--Distortion Theorem.

Let $T\geq T_0$ and consider $\ph_T$ the rescaled map from \eqref{3}. We have to estimate
$$
|\ph_T '(\xi )| 
=\frac{T}{|\ph(T)|}|\ph'(T\xi )|  
\asymp  \frac{T}{|\ph(T)|}|\ph(T\xi )|^{1-\rho}
\quad\text{,}\quad \xi \in Q_1\setminus Q_{1/8}\,.
$$
Assumption \eqref{(4.2)} implies $|\ph(T\xi)| \asymp |\ph(T)|$ which then gives 
$$
|\ph_T '(\xi )| \asymp\frac{T}{|\ph(T)|^\rho}  \succeq 1\quad\text{,}\quad \xi \in Q_1\setminus Q_{1/8}\,,
$$
by \eqref{x.35}. On the other hand, $\frac{T}{|\ph(T)|^\rho} $ is independent of $\xi$. Thus $|\ph_T '(\xi )|\asymp|\ph_T '(1 )|$
and we know from Lemma \ref{14} that $|\ph_T '(1 )|\preceq 1$. Combining all of this gives
$
|\ph_T ' |\asymp  1 $ on $ Q_1\setminus Q_{1/8}$ for $T\geq T_0$. This readily implies that $\b_\infty \equiv 0$.
  \epf 
  
\subsection{H\"older tracts}\label{holder tracts}
  Let $f$ be a model as defined in Definition \ref{tract model} 
or an entire function of class $ \cB$. 
Assume
 that $\Om$ is a single tract of $f$ and that $\ph=\tau^{-1}:\cH\to\Om$ the associated conformal map.
The rectangles $Q_T$ have been introduced in \eqref{2} and $\Om_T=\ph(Q_T)$, $T\geq 1$.
A conformal map $g:Q_1\to U$  is called $(H,\al )$--H\"older if
\beq\label{holder 1'}
|g (z_1)-g(z_2)|\leq H |g' (1)| |z_1-z_2|^\al \quad\text{for all}\quad z_1,z_2\in Q_1\,.
\eeq
The factor $|g'(1)|$ has been introduced in this definition in order to make this H\"older condition scale invariant in the range of $g$. 

\bdfn\label{27}  
The tract $\Tract $ is called H\"older, more precisely $(H,\al )$--H\"older, if \eqref{(4.2)} holds and if there exists $T_0\geq 1$ such that the map $\ph \circ T : Q_1 \to \Tract _T$
satisfies  \eqref{holder 1'} for every $T\geq T_0$. We say that $f$ has H\"older tracts if for some $R\geq 1$
   the components of $f^{-1}(\C\sms\ov \D_R )$ are H\"older tracts.
\edfn

\fr
The main point of this definition is that the tract $\Om$ is exhausted by a family of uniformly H\"older domains $\Om_T$. We shall prove the following.

\blem\label{l1_2017_07_02}
If $f$ is a model or an entire function of class $\cB$ and if $\Om$ is a single H\"older tract of $f$, then
\beq\label{x.15}
|\ph(T)|\asymp \diam(\Om _T) \asymp |(\ph\circ T)'(1)|\quad , \quad T\geq 1\,,
\eeq
where here and in the sequel, by $\ph\circ T$ we mean the map given by the formula
$$
\ph\circ T(z)=\phi(Tz).
$$
\elem

\bpf
The inequality $\diam(\Om _T) \lek |(\ph\circ T)'(1)|$ is immediate from \eqref{holder 1'} while the inequality $\diam(\Om _T) \gek |(\ph\circ T)'(1)|$ is immediate from Koebe's $\frac14$--Distortion Theorem. Hence
$$
\begin{aligned}
|\ph(T)|
&\le |\ph(T)-\ph(1)|+|\ph(1)|
\le H|(\ph\circ T)'(1)|\Big|1-\frac1T\Big|^\a+|\ph(1)| \\
&\lek |\ph(1)|+\diam(\Om_T) \\
&\lek \diam(\Om _T), 
\end{aligned}
$$
where the last inequality was written assuming that $T\ge T_0$ is large enough, say $T\ge T_1\ge T_0\ge 1$. We are thus left to show that 
\beq\label{1mue1}
\diam(\Om_T)\lek |\ph(T)|.
\eeq
Indeed, proving this inequality, it follows from \eqref{holder 1'} and the, already proven, right--hand side of \eqref{x.15}, that
$$
\diam\(\Om_{8^{-q}T}\)\le \frac14\diam(\Om_T)
$$
with some integer $q\ge 1$ and  $T\ge T_1$ large enough, say $T\ge T_2\ge T_1$. Therefore, there exist two points $z_1, z_2\in \Om_T\sms \Om_{8^{-q}T}$ such that
$$
|\ph(Tz_2)-\ph(Tz_1)|\ge \frac14\diam(\Om_T).
$$
So, applying \eqref{(4.2)} we get that
$$
\begin{aligned}
\diam(\Om_T)
&\le 4|\ph(Tz_2)-\ph(Tz_1)|
\le 4|\ph(Tz_2)|+|\ph(Tz_1)| \\
&\le 4(M^q|\ph(T)|+M^q|\ph(T)|)\\
&=8M^q|\ph(T)|,
\end{aligned}
$$
whence formula \eqref{x.15} constituting Lemma~\ref{l1_2017_07_02} is established.
\epf

\brem\label{r1_2017_07_14}
Invoking Lemma~\ref{l1_2017_07_02} we conclude that for a H\"older tract $\Om$, the functions $\ph\circ T$ are uniformly H\"older on $Q_1$ if and only if the rescaled functions $\ph_T$ satisfy uniformly \eqref{holder 1'} without the factor $|g'(1)|$. 
\erem

\fr We also note that if the components of $f^{-1}(\D_{R_0}^*)$ are H\"older for some $R_0\geq 1$ then the components of $f^{-1}(\D_R^*)$ are H\"older for all $R\geq R_0$. 
A very important feature of H\"older tracts is expressed by the following.

\bprop\label{28}
All models or functions $f\in\cB$ with finitely many H\"older tracts have negative spectrum and 
\beq\label{57}
1\leq \tstar \leq \HypDim (f)\leq 2 \; .
\eeq
In addition, if the corresponding H\"older exponent $\a\in (1/2,1]$, then $\tstar<2$. 
\eprop
\bpf
A classical argument (see \cite{PommerenkeBook} or the proof of Proposition~3.3 in
\cite{My-HypDim}) applies word for word showing that 
\beq\label{5_2017_07_07B}
\b_\infty (t+s)\leq (1-\al ) s+ \b_\infty (t),
\eeq
where, we recall, $\a$ is a H\"older exponent of the tract $\Om$. Therefore, 
\beq\label{5_2017_07_07}
b_\infty(t+s)\le b_\infty(t)-\a s.
\eeq
Thus $b_\infty(t)<0$ for all $t>\tstar$ which shows  that $f$ has negative spectrum.

As explained before, the paper  \cite{My-HypDim} also employs the same $b_\infty$--function
along with the zero $\tstar$ and  Theorem 1.1 of \cite{My-HypDim} shows that  \eqref{57}
holds for H\"older tracts.  

The second to last assertion of Proposition~\ref{p1_2017_07_10} is that $b_\infty(0)=1$. So, with $t=0$ and $s=2$, it follows from \eqref{5_2017_07_07} that $b_\infty(2)\le 1-2\a<0$ whenever $\a>1/2$, and thus that $\tstar<2$. The proof is complete.
\epf

All the elementary functions enjoy the property that $\tstar=1$ but in general  H\"older tracts are fractal in the sense that 
$$
\tstar>1\,.
$$ 
Models with fractal tracts have been considered in \cite{My-HypDim}. A particular family of entire functions having fractal tracts is studied in detail in the forthcoming Section \ref{PF functions}.
  
\subsection{Functions of infinite order} \label{31}
Let us finally consider one other family of examples having totally different behavior than the preceding ones.
 They have tracts that are not H\"older, they are of infinite order and also the family of rescalings $\cF_\Om$ has only constant limit functions. Nevertheless, we will see that they have negative spectrum and thus they are first examples of infinite order for which the thermodynamic formalism is developed. 

 Consider functions $f\in \cD$, no matter whether entire or model,
 having the following properties:
 
 \begin{itemize}
 \item[-] $f$ has negative spectrum.
 
 \sp\item[-] $f$ has a H\"older tract $\Om_f\subset \{\Re z\geq 3\} \subset \cH$.

 \end{itemize}
We will associate to such a function  a model function $F$ defined on $\Om_F=\log (\Om_f)$, $\log$ meaning any, or even finitely many, arbitrary branches of the logarithm. The definition of $F$ is this.
$$
F:= f\circ\exp:\Om_F\to \C\,.
$$
To such a function Theorem~\ref{thm intro} applies since we have the following.
 
\bprop\label{z.8}
The infinite order function $F= f\circ \exp : \Om_F\to \C$ belongs to $\cD$  and has negative spectrum with $\tstarF = \tstar$.
\eprop
 
\bpf 
The disjoint type property follows since $F^{-1}(\D^*)=\Om_F\subset \D^*$. Let $\ph = \tau^{-1} : \cH\to \Om_f$ be conformal such that $f=e^{\tau}$ on $\Om_f$. Then 
$$
 F=\exp \circ (\tau \circ \exp). 
$$
 It suffices to consider the case where $\Om_F$ is a single tract so that $ \tau \circ \exp: \Om_F \to \cH$ is a conformal map
 with inverse $ \Phi = \log \circ \,\ph$.
Since $f\in \cD$ it satisfies  \eqref{(4.2)} and since $\Om_f \subset \{\Re z\geq 3\}$ we have
$$
\frac{|\Phi(\xi_1)|}{|\Phi(\xi_2)|}
=\frac{\log |\ph (\xi _1)|}{\log |\ph (\xi _2)|} \leq 1+\frac{\log M}{\log |\ph (\xi _2)|}
 \leq 1+ \frac{\log M}{\log 3}
  \quad \text{for all} \quad \xi_1,\xi_2\in Q_T\setminus Q_{T/8}\,.
$$
Thus $\Phi$ satisfies \eqref{(4.2)}, completing the argument that $F$ is in $\cD$.
  
  It remains to estimate $\b_{\infty , F}$. For $T\geq 1$, $\Phi_T = \frac{1}{|\Phi (T)|} \log \circ \ph \circ T $ hence
  $$
 \Phi_T'(\xi ) = \frac{T}{|\Phi (T)|} \frac{\ph'(T\xi )}{\ph (T\xi )} = \frac{1}{|\log \ph (T)|}\frac{|\ph (T)|}{\ph (T\xi)}\, \ph_T'(\xi )
  $$
 thus
 $$
 | \Phi_T'(\xi )  | \asymp \frac{1}{\log |\ph (T)|}| \ph_T'(\xi )|\quad \text{for every} \quad \xi \in Q_1\setminus Q_{1/8}\,.
 $$
 The factor $|\ph (T)|$ can be estimated as follows. Still since $\Om_f \subset \{\Re z\geq 3\}$ we have 
 $|\ph(T)|\geq 3$. On the other hand we have from \eqref{bdd 3}
 $$
 |\ph(T)| \lek |\ph(T)-\ph(1)| \leq T^2 |\ph'(1)|\,.
 $$
 It follows that there exists a constant $C\geq 0$ such that
 $$
 \frac{e^{-C}}{\log T}\int _{I} |\ph_T'(r+iy)|^t dy\leq \int _{I} |\Phi_T'(r+iy)|^t dy\leq e^C \int _{I} |\ph_T'(r+iy)|^t dy
 $$
 for every $r\in (0,1)$ and $T\geq 1$. Setting now $T=1/r$, taking logarithm of the preceding inequality and dividing it then by $\log 1/r$  shows that 
$$ \frac{-C-\log \log 1 /r}{\log 1/r} + \b _{\ph _{1 /r}}(r,t)\leq
 \b_{\Phi_{1/r}}(r,t)\leq  \b _{\ph _{1 /r}}(r,t).
 $$
 Since the first term on the left hand side of this inequality goes to zero as $r\to 0$ we get
   $ \b_{\infty , F}(t)=  \b_{\infty , f}(t)$ from which the Proposition follows.
\epf

 
\section{Quasiconformal invariance of H\"older tracts}\label{QIHT}
Quasiconformal maps have good H\"older continuity properties and thus preserve H\"older tracts. Let us make this precise.

\blem\label{w.4}
Let $g, f\in \cB$ have finitely many tracts and assume that there is $R\geq 1$ such that all the connected components of $g^{-1}(\D_R^*)$
are H\"older. If $\Phi:\C\to\C$ is a quasiconformal homeomorphism such that 
$$
f^{-1}(\D_R^*) = \Phi \left(g^{-1}(\D_R^*) \right)
\quad \text{and}\quad 
f\circ \Phi = g \quad \text{on}\quad g^{-1}(\D_R^*)
$$
then all the connected components of $f^{-1}(\D_R^*)$ are H\"older.
\elem

\bpf
It suffices to consider the case where the functions have just one tract
$\Om_g= g^{-1}(D_R^*)$ and $\Om_f =\Phi (\Om_g)$. We may also assume without loss of generality that $R=1$. Then there are conformal maps $\ph_g:\cH\to\Om_g$ and $\ph _f:\cH\to\Om_f$ such that, with appropriate holomorphic branches of logarithms,
the holomorphic maps $\log g:\Om_g\to\cH$ and $\log f:\Om_f\to\cH$ are the respective inverses of  $\ph_g$ and  $\ph _f$, and, in addition,
$$ 
\ph _f = \Phi \circ \ph _g\,.
$$
By our hypotheses, $\ph_g$ satisfies the conditions of Definition \ref{27} and 
we have to show that $\ph_f$ does it too.
The condition \eqref{(4.2)} is satisfied by $\ph_g$ and $|\ph_g (T)|\to \infty$ as $T\to\infty$. Thus
\beq\label{5_2017_07_10}
|\ph_g(\xi_1) -\ph_g(0)|\asymp  |\ph_g(\xi_1)|\leq M  |\ph_g(\xi_2)|\asymp  |\ph_g(\xi_2) -\ph_g(0)|
\eeq
for all $ \xi_1,\xi_2\in Q_T\setminus Q_{T/8}$, $T\geq 8$. Since the map $\Phi$ is quasiconformal, it is quasisymmetric. This means that there exists a homeomorphism $\eta: [0,\infty )\to[0, \infty)$ such that $|a-b|\leq t |a-c|$ yields 
$$
|\Phi (a)-\Phi (b)|\leq \eta (t)|\Phi (a)-\Phi (c)|.
$$
Along with \eqref{5_2017_07_10}, this gives that 
$$ 
\begin{aligned}
|\ph_f(\xi_1)|
&=|\Phi( \ph_g (\xi _1 ))|
\comp |\Phi( \ph_g (\xi _1 ))-\Phi(\ph_g (0))| \\
&\leq \eta(M') |\Phi( \ph_g (\xi _2 ))-\Phi(\ph_g(0))| \\
&\comp|\Phi( \ph_g (\xi _2 ))| \\
&=|\ph_f(\xi_2)|
\end{aligned}
$$
for all  $\xi_1, \xi_2\in Q_T\setminus Q_{T/8}$, $T\geq 8$, where $M'$ is a constant witnessing the comparability of the very left and the very right sides of  \eqref{5_2017_07_10}.
In other words, $\ph_f$ satisfies \eqref{(4.2)}.

We know that for some $T_0\geq 1$ the family of rescalings $\ph_{g,T}$, $T\geq T_0$, of $\ph_g$ is uniformly H\"older and it remains to show that
$$
\ph_{f,T} = \frac{1}{|\Phi (\ph_g (T))|}\,  \Phi \circ |\ph_g (T)| \circ \ph _{g,T}
  \quad , \quad       T\geq T_0\, , 
$$
  has
the same property. All the mappings 
$$
\hat g_T:= \frac{1}{|\Phi (\ph_g (T))|}\,  \Phi \circ |\ph_g (T)| 
$$ 
are $K$--quasiconformal, where $K$ is the quasiconformal constant of $\Phi$ and they are normalized by $\hat g_T(\infty)=\infty$. We shall prove the following.

\sp{\bf Claim:} There exists a constant $\ka\in(0,1]$ such that 
$$
|\hat g_T(0)|\leq \ka \leq 2\ka \leq |\hat g_T(1)| \leq 1/\kappa \quad \text{ for all $T\geq T_0$.}
$$

\bpf
We have 
$$
|\Phi(|\ph_g (T)|)|\comp |\Phi(|\ph_g (T)|)-\Phi(0)||
\ \  \  {\rm and } \  \  \
|\Phi(\ph_g (T))|\comp |\Phi(\ph_g (T))-\Phi(0)||,
$$
and obviously
$
|\ph_g (T)-0|=\big||\ph_g (T)|-0\big|.
$
Therefore, invoking again quasisymmetricity of $\Phi$, witnessed by the  homeomorphism $\eta: [0,\infty )\to[0, \infty)$, we consecutively get
$$
\frac1{\eta(1)}
\le \frac{|\Phi(|\ph_g (T)|)-\Phi(0)||}{|\Phi(\ph_g (T))-\Phi(0)||}
\le\eta(1)
$$
and 
$$
\frac{|\Phi(|\ph_g (T)|)|}{|\Phi(\ph_g (T))|} \asymp 1.
$$
This means that
$
 |\hat g_T(1)|\asymp 1.
$
The proof of the Claim is complete.
\epf
In conclusion, 
$$
\cG:=\{g_T: T\geq T_0\}
$$ 
is a uniformly quasiconformal and normalized family. By Remark~\ref{r1_2017_07_14} there exists $R>1$ such that 
$$
\ph_{g,T}(Q_1)\subset \ov\D_R
$$
for all $T\geq T_0$. These two facts imply (see Theorem 4.3 in \cite{LV-1973}) that the family $\cG$ restricted to $\ov\D_R$ is uniformly H\"older.
Therefore, $\ph_{f,T}$ is  uniformly H\"older as a composition of two  H\"older functions whose H\"older exponents and constants do not depend on $T\geq T_0$.
\epf

We provide two important applications  of the quasiconformal invariance of H\"older tracts. The first one we present right now and it concerns  quasiconformal approximation. The second application will be in Section \ref{section S} on analytic families of functions in  the Speiser class $\cS$.

\subsection{Quasiconformal approximation} As already mentioned, Bishop 
  \cite{Bishop-EL-2015, Bishop-S-2016} considered quasiconformal approximations of most general models where $\Tract$ can be an arbitrary union of simply connected unbounded domains. Keeping our definition of a model, the following result is a simplified version of Theorem~1.1 in \cite{Bishop-EL-2015}.

\bthm[\cite{Bishop-EL-2015}]\label{bishop}
Let $(\tau , \Tract )$ be a tract model and $f=e^\tau$ the corresponding function. Fix $R>0 $. Then there exist an entire function $F\in \cB$
and a quasiconformal map $\psi :\C\to\C$ with $\psi$ conformal outside of $ \{z\in\C:R <\Re \tau (z) < 2R\}$ such that 
$$
e^\tau = F\circ \psi
$$ 
on $ \Om(2R)=\{z\in\C:\Re \tau (z) >2R\}$. Moreover, the components of $\{|F|>e^R\}$ are in a $1$-to-$1$ correspondence with the components of $\Om$ via $\psi$.
\ethm

If the initial model is of disjoint type one can adjust $R>0$ such that 
\beq\label{rho} 
\J_f\subset  \Om(3R)=\{z\in\C:\Re \tau (z) >3R\}\,.
\eeq
Then we can assume that $F$ also is of disjoint type 
since otherwise it suffices to compose $\psi$ with an affine map. So, we can consider for the map $F$ the set of tracts $\Om_F=\psi (\Om (2R))$ and suppose that $\J_F\subset \psi (\Om ( 3R))$.

\bprop \label{z.10}
Suppose $f\in \cD$ is a model having only H\"older tracts and suppose that $F$ is a disjoint type entire function given by Bishop's Theorem \ref{bishop} with $R>0$ small enough such that \eqref{rho} holds. Then, $F\in \cD$, the tracts $\Om_F$ of $F$ are H\"older and consequently $F$ has negative spectrum.
\eprop

\bpf
Follows directly from Lemma \ref{w.4}.
\epf

\smallskip

\section{Poincar\'e functions}\label{PF functions}

In this section, we consider a, quite particular, family of entire functions.
 They are obtained by linearization of a polynomial 
 at a repelling fixed point and are often called linearizers or Poincar\'e functions (see \cite{DO08, MP12, ER15}).
 We first consider general linearizers and then show that these functions have negative spectrum if and only if the corresponding 
 polynomial is topologically Collet-Eckmann (TCE).
 
 \subsection{Linearizers: the general case.}
 Let $p:\oc\to\oc$ be a polynomial having a repelling fixed point $z_0$ with multiplier $\l = p'(z_0)$.
By the Koenigs-Poincar\'e linearization theorem there exists an entire function $f:\C\to\C$ such that
\beq\label{x.3}
  f(0)=z_0\quad \text{and} \quad f(\l z) = p\circ f(z) \quad , \quad z\in\C.
\eeq
We then call $f$ linearizer of $p$ or, more precisely, linearizer of $p$ at the fixed point $z_0$.

\smallskip

\brem
One could consider here a much more general family. 
Instead of linearizing the dynamics at a repelling fixed point one can consider
limits of rescalings at conical limit points. For example, if $p$ is a hyperbolic 
polynomial with connected Julia set then there exist entire linearizers at any point of the Julia set 
$\J_p$ (see \cite[Theorem 2.10]{BFU02}).
\erem

We only will consider the case where the Julia set of $p$ is connected. Then $\P(p)$, the post--singular set of $p$,
 is a subset of the filled Julia set and thus a subset of the complement of $A_p(\infty )$, the attracting bassin of infinity.
Since the  set of singular values of $f$ is equal to the post--singular set of $p$ (see Proposition~3.2 in \cite{MP12})
it follows that  $f$ belongs to class $\cB$. Up to normalization we can assume that $\J_p\subset \D$ or, equivalently,
that $\D^* \subset A_p(\infty )$. Then $S(f)\subset \D$ and we can consider the set of tracts $\Om =f^{-1} (\D^*)$.
Conjugating eventually $p$ by an affine map,
we also may assume  without loss of generality that 
$$\text{$z_0\neq 0\; $ and that $\;0\not\in A_p(\infty)$.}$$

Linearizers of polynomials are entire functions of finite order \cite{Valiron1913}
 and the Denjoy--Carleman--Ahlfors Theorem asserts that finite order functions have only finitely many tracts
 $\Om_j \subset \Om$, $j=1,...,N$. Let $\hat \Om_j$ be the connected component of $ f^{-1} (A_p(\infty ))$
 that contains $\Om_j$.
 It follows from $p$--invariance of $A_p(\infty )$ along with \eqref{x.3} that multiplication by $\l$ permutes 
 the finitely many components $\hat \Om_j$ (see also Proposition 2.1 of \cite{Hub1993} for a different approach to this fact).
 Replacing $p$, hence $\l$, by an iterate we may assume 
 that  
 $$\text{$\l \hat\Om_j =\hat \Om_j \; $ for every $\; j=1,...,N$.}$$
 It  thus suffices to consider in the following a single tract. We will call it $\Om$, thus from now on $\Om$ will be a connected component of  $ f^{-1}(\D^*)$ and $\hat \Om $ will be the component of 
$f^{-1}(A_p(\infty ))$ that contains $\Om$. Recalling that $0\notin A_p(\infty)$, we see that the function $f$ restricted to $\hat \Om$ is again of the form 
$
f=e^{\ph ^{-1}},
$
where $\ph: \hat \cH \to \hat \Om$ is a conformal homeomorphism and $\Om=\ph(\cH)$.
 
 \smallskip
 
Since $\J_p$ is connected, there exists a
conformal map $h:\D^*\to A_p(\infty)$ such that 
 \beq\label{x.9}
 h(1)=z_0 \quad \text{and}\quad
 p\circ h (z) = h(z^d) \quad , \quad z\in\D^*,
\eeq
where $d\ge 2$ is the degree of the polynomial $p$. 
Notice that here $h(1)$ must be properly defined since the Julia set
$\J_p$ may fail to be locally connected and then $h$ has no continuous extension to the boundary. 
However, since $z_0$ is repelling there exists 
an external ray $\cR_t=h(\{re^{it}\, ; \; r>1\})$ who lands at $z_0$. Such a ray can be obtained
by constructing first a $\l$--invariant curve in the tract $\hat \Om$ which then can be projected
down by $f$ (see Section 2.3 of \cite{Hub1993}). 
In our situation we may assume that $t=0$
and then the normalization in  \eqref{x.9} is defined with the radial limit $\lim _{r\downto 1}h(r)=z_0$.

The map $h$ can be lifted, via the exponential map, to a  conformal homeomorphism 
$$
H:\cH \to \hat \cH:=\exp^{-1} (A_p(\infty ))
$$ 
that commutes with translation by $2\pi i$, i.e. 
$
H(z+2\pi i) = H(z)+2\pi i$ for all $ z\in \cH$,
and such that
\beq\label{x.13}
\ph\circ H (0)=0 
\quad \text{and}\quad 
\exp\circ H=h\circ\exp\quad on \quad \cH\, .
\eeq

\blem\label{cl1_2017_07_19}
The inverse conformal map $H^{-1}:\hat \cH\to\cH$ is bi--Lipschitz on the half--space $\{\Re z >s\}$, whenever $s\in \R$ is such that $\{z\in\C:\Re z\ge s\} \subset \hat \cH$. 
\elem

\bpf
Indeed, with $\tilde h (z):= 1/h(1/z)$, $z\in \D$,
$$
H'(z) = \frac{e^{-z}}{\tilde h (e^{-z})} \tilde h '(e^{-z}) \longrightarrow 1 \quad \text{when}\quad \Re z \to +\infty\,.
$$
This and the $2\pi i$--periodicity of $H$ imply the announced bi--Lipschitz property.
\epf

\fr Since $\ov\D^*\subset A_p(\infty)$, there exists $t<0$ such that 
\beq\label{2_2017_07_19}
\ov \cH \subset \{\Re z >t\}\subset \hat \cH\,.
\eeq
Consequently, $H^{-1}$ is uniformly bi--Lipschitz, say with constant $L\geq 1$, on all the rectangles $Q_T$, $T\geq 1$.

\smallskip

\begin{figure}[h]
   \includegraphics[height=7.7cm]{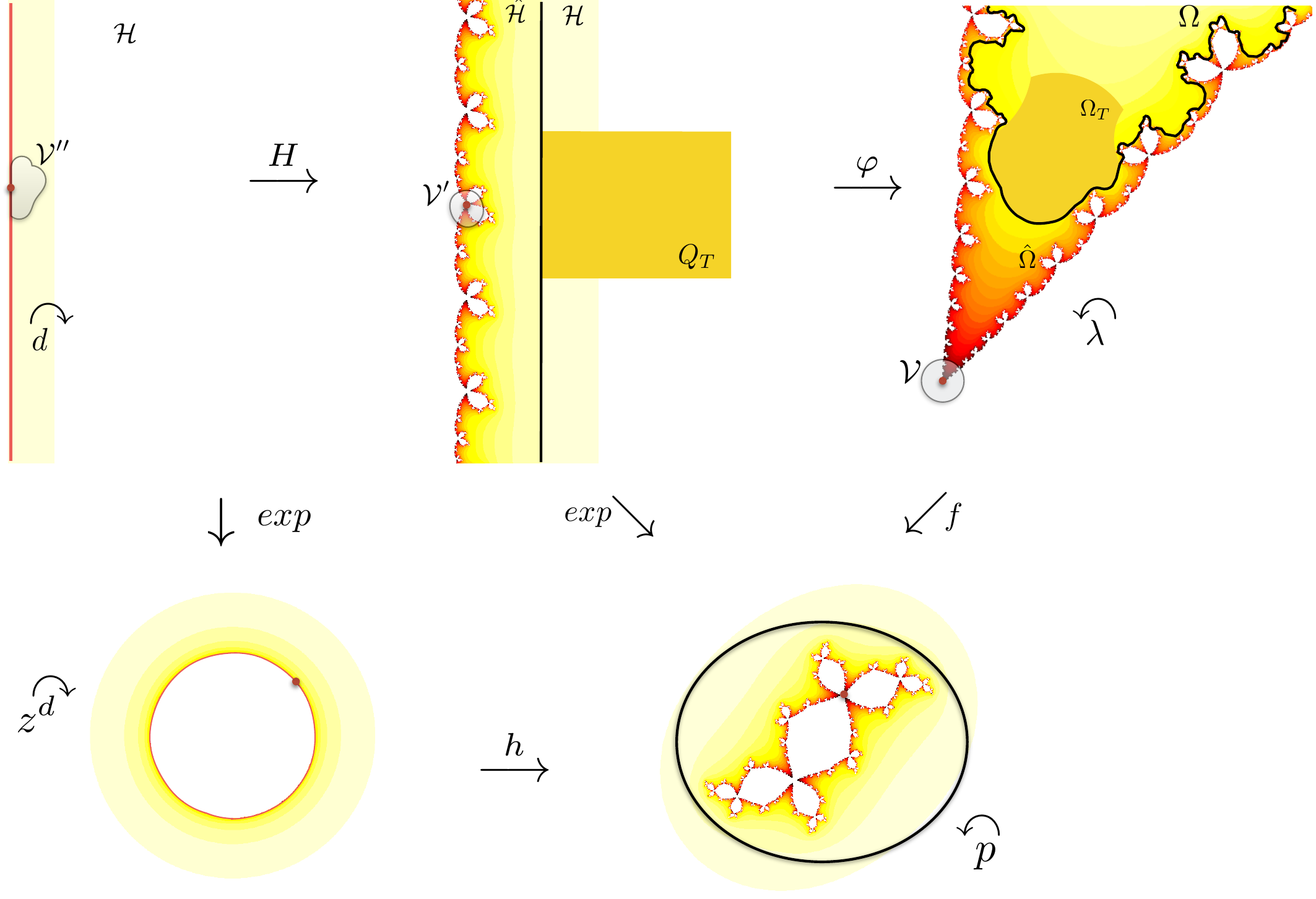}
      \caption{Linearizing Douady's Rabbit}
      \label{Figure 4}
\end{figure}

Recall that $\l$ is the multiplier of $p$ at the repelling fixed point $z_0$. Since $f$ conjugates multiplication by $\l$ and the polynomial $p$, since $h$ satisfies \eqref{x.9} and the right hand sided part of \eqref{x.13}, and since the exponential map lifts $z\mapsto z^d$ to multiplication by $d$, we get
$$
f\big(\ph\circ H (dz)\big) =h\big(e^{dz}\big)=p\big(h(e^z)\big)=
p\big(f(\ph\circ H)(z)\big)= f\big(\l (\ph\circ H)(z)\big) \,
\text{  for all $z\in \cH$.} 
$$
Along with the left hand side of \eqref{x.13} this implies that
\beq\label{x.17}
 \big( \ph\circ H \big)(d z ) = \l\,  \big(\ph \circ H\big) (z ) \quad , \quad z \in \cH\,.
 \eeq
 Indeed, both $ \ph\circ H\circ d$ and $\l(\ph\circ H)$ are conformal mappings from $\cH$ onto $\hat \Om$ 
 and \eqref{x.17} holds near the origin. In conclusion, we have the commutative diagram of Figure \ref{Figure 4}. 
 
As a final preliminary remark we recall that a linearizer is essentially unique. Indeed, if
 $f$ has the property \eqref{x.3}, then every other solution of \eqref{x.3} is of the form
$f_\ka=f\circ \kappa$ 
for some $\ka \in \C^*$.

\bthm\label{x.6 general}
If $p$ is a polynomial with connected Julia set, if $f_\ka=f\circ \ka$ is a linearizer of $p$ at a repelling fixed point $z_0$
 and if $\ka$ is sufficiently small, then, up to normalization, 
$\cJ_{f_\ka}\subset \D^*\subset A_p(\infty)$ and the following holds:
\ben
\item $f_\ka \in\cD$.

\item $\Th _{f_\ka}=\HypDim(p)$, the hyperbolic dimension of $p$.

\item The tracts of $f_\ka$ are fractal in the sense that 
$$\Th _{f_\ka} >1$$
if and only if $p$ is not an exceptional polynomial.
\een
\ethm

Recall that a polynomial is called exceptional if it is either of the form $z\mapsto z^d$, $d\ge 2$, or it is a Tchebychev polynomial. Such polynomials are special in the sense that Zdunik \cite{Zdunik90} has shown that they are the only polynomials for which the harmonic measure, viewed from infinity, is not singular with respect to the natural Hausdorff measure of the Julia set.

   \smallskip

Theorem \ref{x.6 general}  will be shown in two steps. It will be a consequence of the next Lemma along with Proposition \ref{x.36'}.

\blem \label{x.36}
If $p$ is a polynomial with connected Julia set and if $f_\ka=f\circ \ka$ is a linearizer of $p$ at a repelling fixed point $z_0$ then 
$f_\ka \in\cD$ provided  $\ka$ is sufficiently small 
and, up to normalization, 
$\cJ_{f_\ka}\subset \D^*\subset A_p(\infty)$.
\elem

\bpf 
By the above discussion, we already know the following facts: the linearizers $f_\ka$ are of finite order and thus 
have only finitely many tracts, they belong to class $\cB$ and, normalizing $p$ if necessary,  
$\cJ_{f_\ka}\subset \D^*\subset A_p(\infty)$.
Furthermore, for small values of $\ka$ the function $f_\ka$ is of disjoint type and, in order to verify that it belongs to class $\cD$,
it remains to verify \eqref{(4.2)}. Notice that this condition is equivalent to the following: there exists $C>1$ and $M<\infty$ such that
\beq \label{(4.2')}   
\text{$\forall \xi , \xi '\in \cH \; $ with }\;  \frac 1C < \frac{|\xi|}{|\xi'|}<C \quad \text{have} \quad
|\ph (\xi )|\leq M |\ph (\xi ' )|\,.
\eeq
Moreover, if \eqref{(4.2')}  holds for some $C>1$ then it does hold for all $C>1$ (and for some $M=M(C)<\infty$).
Since here $\ph = (\ph \circ H)\circ H^{-1}_{|\cH}$ it suffices to verify separately \eqref{(4.2')} for $H^{-1}_{|\cH}$ and for $\ph \circ \H$.
Then clearly the composition $ (\ph \circ H)\circ H^{-1}_{|\cH}$ also satisfies \eqref{(4.2')} (with different constants).

\smallskip

For $H^{-1}_{|\cH}$ this directly follows from the bi-Lipschitz property in Lemma \ref{cl1_2017_07_19}: let $C>1$ and 
$\xi , \xi '\in \cH$ with $\frac 1C < \frac{|\xi|}{|\xi'|}<C$. Then, for all $\xi \in \cH$, 
\beq\label{new 11}
\frac 1L |\xi| - |H^{-1}(0)|\leq |H^{-1}(\xi)|= |H^{-1}(\xi) -H^{-1}(0)+H^{-1}(0)|\leq L |\xi | +|H^{-1}(0)|
\eeq
hence,  if $|\xi |, |\xi'| \geq 2L|H^{-1}(0)|$,
$$
\frac 1C \frac{\frac 1{2L}}{L+\frac 1{2L}}\leq \frac{|H^{-1}  (\xi)|}{|H^{-1} (\xi')| }\leq 
\frac{\left( L+\frac 1{2L}\right) |\xi|}{\left( \frac 1L-\frac 1{2L}\right) |\xi'|} \leq \frac{\left( L+\frac 1{2L}\right) }{\left( \frac 1L-\frac 1{2L}\right)}C
$$
which implies the required inequality with $M=\frac{L+1/(2L)}{1/L - 1/(2L)}C$ in this case.
If one of the points $\xi , \xi '$ has modulus smaller than $2L|H^{-1}(0)|$ then both have modulus smaller than $s=2L|H^{-1}(0)|C$
and then the required inequality follows since for the compact subset $K=\ov{\D_s\cap \cH}$ of $\hat \cH$
one clearly has $\inf _{K}|H^{-1} | >0 $ and $\sup _{K}|H^{-1} | <\infty$.

\smallskip

It remains to proof \eqref{(4.2')}  for the map $\ph \circ H$. This map conjugates multiplication by $d$ and multiplication by $\l$
(see \eqref{x.17}). The set $A=\{ \xi \in \cH \; , \; 1< |\xi | \leq d\}$ is a fundamental set for the multiplication by $d$ in $\cH$,
$A'= \ph \circ H (A)\subset \hat \Om \subset \C^*$ and $A'/<\l >$ is an annulus of the torus $\C^*/<\l >$
(whose precise description is given in \cite{Hub1993}). Consequently, $\ov {A'}$ is a compact subset of $\C^*$
as well as the set $\l^k \ov {A'} = \ov{\ph \circ H (d^k A)}$, $k\in \Z$.

Let $C>1$ and let $\cA=\Big\{\frac 1C \leq |\xi | \leq Cd\; , \; \xi \in \cH\Big\}$. This set can be covered by finitely many sets $d^k A$
and thus $\ov{\ph \circ H  (\cA)}$ is also compact in $\C^*$. Therefore,
there exists $0<\underline{m}\leq \overline{m}<\infty$ such that
\beq\label{new 1}
\underline{m} \leq |\ph \circ H( \xi )| \leq \overline{m} \quad \text{for all}\quad \xi \in \cA=\Big\{\frac 1C \leq |\xi | \leq Cd\; , \; \xi \in \cH\Big\}.
\eeq 
Now, let $\xi, \xi'$ be like in \eqref{(4.2')} and choose $k\in \Z$ such that $1<|d^k \xi |\leq d$.
Then $d^k\xi, d^k\xi '\in \cA$ and thus \eqref{x.17} and \eqref{new 1} imply that
$$
\frac{|\ph\circ H (\xi) |}{|\ph\circ H(\xi') |} = \frac{|\ph\circ H (d^k\xi) |}{|\ph\circ H(d^k\xi') |} \leq \frac{\ov m}{\underline{m}}.
$$
This completes the proof since it shows that  $\ph \circ H$ also satisfies \eqref{(4.2')}.
\epf

\smallskip

The second and final step of the proof of Theorem \ref{x.6 general} is the following.

\bprop \label{x.36'}
Let $p$ be a polynomial with connected Julia set and let $f$ be a linearizer of $p$ that satisfies the conclusions of Lemma \ref{x.36}.
Then,
\ben
\item $\Th _{f}=\HypDim(p)$, the hyperbolic dimension of $p$.
\item The tracts of $f$ are fractal in the sense that 
$$\Th _{f} >1$$
if and only if $p$ is not an exceptional polynomial.
\een
\eprop

A key point in the proof of this result is that we will be able to relate the $\b_\infty$ function of the rescalings to the classical integral means spectrum
$$
\b_h(t)=\limsup_{r\to1^+}\frac{\log \int _{|z|=1}|h'(rz)|^t|dz|}{-\log(r-1)}
$$
of the Riemann map $h:\{|z|>1\}\to A_p(\infty)$ of \eqref{x.9}. For this function there is a formula holding for all polynomials with connected Julia sets (see \cite{BMS03}, see also \cite{PUbook} for the expanding case):
 \beq \label{x.21}
 \b_h(t)-t+1 = \frac{\P(t)}{\log d}
 \eeq
 where $d:=\deg(p)$ and $\P(t)$ is the topological pressure of the potential $-t\log|p'|$ with respect to the polynomial $p$. In fact $\P(t)$ is the tree pressure in  the general non-expanding case, see \cite{PRZ99,PRS04}. Since  this formula does not directly hold here, we provide its suitable variant along with all the details in the Appendix \ref{y.1}.

\bpf In order to relate these spectra we have to study in detail the rescaled functions $\ph_T$ and thus $\ph$.
We start with a preliminary observation. If $\exp_*^{-1}$ is the holomorphic inverse branch of the exponential map, defined near $z_0$, such that $\exp_*^{-1}(z_0)=\tau(0)$,
then $\exp_*^{-1}\circ f$ extends $\tau = \ph ^{-1}$ to some bounded open neighborhood, call it $V$, of the origin.
Since $f'(0)\neq 0$, the neighborhood $V$ can be chosen such that
$\tau: V\to V'=\tau(V)$ is bi--Lipschitz. Then
\beq\label{new 14}
|\ph'| \asymp 1 \quad \text{ in }\quad V'.
\eeq
Again \eqref{x.17} shows that
$
\ph = (\ph\circ H)\circ H^{-1}_{|\cH} = \l ^N (\ph\circ H)\circ d^{-N}\circ H^{-1}  $ and thus 
$$
\ph \circ T _{|Q_1} =  \l ^N (\ph\circ H)\circ (d^{-N}\circ H^{-1}\circ T)_{|Q_1} : Q_1\to \Om_T
$$
for every $N\in \Z$.  Given $T\geq 1$, we fix this integer as follows. Let $A$ be again the fundamental set for the multiplication by $d$ in $\cH$
of the proof of Lemma \ref{x.36} and $A'= \ph\circ H (A)$. Let then $N_0=N_0(T)\in \Z$ be the unique integer such that
\beq\label{new 31}
\l^{-N_0} \ph (T)\in A'\quad \text{or, equivalently, such that} \quad d^{-N_0} H^{-1} (T) \in A.
\eeq
Since $H^{-1}_{|\cH}$ is bi-Lipschitz we have, by the same argument than in \eqref{new 11}, that
\beq\label{new 12} T\asymp d^{N_0} \quad \text{uniformly in $T\geq1$}\eeq
and thus we see that the map $d^{-N_0}\circ H^{-1}\circ T $ is uniformly bi-Lipschitz on $Q_1$.
On the other hand, \eqref{new 1} implies that  
\beq\label{new 13}
\underline{m} \leq |\l |  ^{-N_0} |\ph (T)| \leq \overline{m} 
\eeq
thus, using  \eqref{(4.2)}, 
$$
\frac{\underline{m}}{M}\leq |\l |^{-N_0}|\ph (\xi )|\leq M\overline{m} \quad \text{for all}\quad \xi \in Q_T\setminus Q_{T/8}
\; , \; T\geq 1.
$$
Therefore, there exists an integer $n_0$ such that for every $T\geq 1$ 
$$
\l^{-n_0-N_0} \ph (Q_T\setminus Q_{T/8} ) \subset V.
$$
Let $N=n_0+N_0$. Then, $n_0$ being fixed, \eqref{new 12} and \eqref{new 13} still hold (with different constants) if we replace $N_0$ by $N$ and the maps $G_T= d^{-N}\circ H^{-1}\circ T $ are still uniformly bi-Lipschitz.

Let us now consider the rescaled maps
\beq\label{new 30}
\ph_T 
= \frac{1}{|\ph(T)|} \, \ph\circ T 
=\frac{\l^N}{|\ph(T)|}\,  \ph \circ H \circ G_T :Q_1\to  \C \,.
\eeq
We know that $|G'_T|\asymp 1$ on $Q_1$ and that $|\l|^N\asymp |\ph(T)|$. Moreover, $\ph$ is evaluated on the set 
$ Q_{1,T}=H\circ G_T(Q_1)$ and, by the choice of $n_0$ hence of $N$, $ Q_{1,T}\subset V'$.
 But on $V'$ we have \eqref{new 14}
and thus 
\beq \label{z.12}
|\ph_T '| \asymp | H' \circ  G_T | \quad \text{on} \quad Q_1\,.
\eeq
We are now ready to investigate the integral means.
Remember that in the definition of $\b_{\ph_T}$ one integrates over the set 
$
I_r =\big \{r+iy: \; 1 < |y| <2\big\}$, $0<r<1$.
Fix $r\in(0,1)$ and focus on the set
$$
I_r^+=\{z\in I_r: \; \Im (z) >0\}
$$
and, as in the definition of $\b_\infty (r,t)$,  consider in what follows $T\geq \ga / r$. We claim that then there exists $K\geq L\geq 1$,
$L$ the bi-Lipschitz constant of Lemma~\ref{cl1_2017_07_19},  such that
\beq\label{x.34}   
G_T(I_r^+) \subset \Big\{z\in \C: 0<\Re z < Kr\Big\}\quad \text{for every} \quad T\geq \ga/r \,.
\eeq
Indeed, since $H^{-1}(z+2\pi i)=H^{-1}(z)+2\pi i$, we have that $\Re(H^{-1}(z+2\pi i))=\Re(H^{-1}(z))$, and therefore
there exists $C>0$ such that $\Re(H^{-1}(iy))\leq C$ for all $y\in \R$. By Lemma~\ref{cl1_2017_07_19} we have
$$
|H^{-1} (Tr+iy ) -H^{-1} (iy)|\leq  L Tr 
$$
for all $y\in \R$. Thus, for all $z=r+iy\in I_r^+$,
$$
\begin{aligned}
0&<\Re G_T( z) 
=d^{-N}  \Re H^{-1} (Tr+iTy)\\
&\leq d^{-N} \Big| \Re H^{-1} (Tr+iTy)-\Re H^{-1} (iTy)\Big|+d^{-N} \Big|\Re H^{-1} (iTy) \Big| \leq d^{-N} \Big( LTr+C\Big)
\end{aligned}
$$
and since $T\asymp d^N$, by \eqref{new 12}, we get
$$
0<\Re G_T( z )  \preceq Lr+ C/T \leq (L+C/\ga)r \quad , \quad z=r+iy\in I_r^+,
$$
which shows \eqref{x.34}.
Let now $\sg  \subset (\{\Re z = Lr\} )$ be a sufficiently long compact line segment so that $\ga :=G_T^{-1}(\sg )$ is a cross-cut of $ \big\{\xi\in\C: 1\leq \Im \xi \leq 2\big\}$.
For every $k\in \{0,...,[1/r]\}$ set 
$$
a_k(r):= r+i(1 +kr)
$$ 
and let $b_k(r)\in \ga$ with 
$$
\Im (b_k (r)) = 1 +kr.
$$
We can choose the points $b_k$ such that $\Im c_{k+1}> \Im c_{k} $, where $c_k:=G_T(b_k)$. It follows from \eqref{x.34} and Lemma \ref{cl1_2017_07_19}
  that the Hausdorff distance between $I_r^+$ and $\g$ is bounded above by a multiple of $r$ and, in addition, these two sets are disjoint and (recalling that $G_T$ is orientation preserving since conformal) $\Re (w)>r$ for all $w\in\g$. Hence, there exists $\ka>1$ and for every $r\in(0,1)$ and every $k\in \{0,...,[1/r]\}$ there exists a rectangle $\De_k(r)$ 
whose ratio of the longer to the lower edge is uniformly bounded above
such that $a_k(r), b_k(r)\in \De_k(r)$ and $\ka\De_k(r)\sbt \cH$. Therefore we can apply Koebe's Distortion Theorem to the map $\ph_T|_{\ka\De_k(r)}$ to conclude that
$$
|\ph_T '(a_k(r))|\comp|\ph_T '(b_k(r))|
$$
with a comparability constant independent of $r$ and $k$. Then,
$$
\int _{I_r^+}|\ph_T '(\xi )| ^t |d\xi |
\asymp \sum_{k=0}^{[1/r]} |\ph_T '(a_k )| ^t r\asymp
 \sum_{k=0}^{[1/r]} |\ph_T '(b_k )| ^t r
 \asymp
 \sum_{k=0}^{[1/r]} |H '(c_k )| ^t r
$$
where the last comparability sign follows directly from \eqref{z.12}. On the other hand,
$$
\int _\sg |H'(z)|^t |dz| 
\asymp  \sum_{k=0}^{[1/r]}  |H '(c_k )| ^t |c_{k+1} -c_{k}| 
\asymp \sum_{k=0}^{[1/r]}  |H '(c_k )| ^t \, r
$$
since $|c_{k+1} -c_{k}| \asymp|b_{k+1} -b_{k}| \asymp r$. This shows that 
$$
\int _{I_r^+}|\ph_T '(\xi )| ^t |d\xi |\asymp\int _\sg |H'(z)|^t |dz|.
$$
Having this, an elementary calculation (Chain Rule), based on \eqref{x.13} and on the fact that $z_0\neq 0$, yields
\beq\label{40}
\int _{I_r^+}|\ph_T '(\xi )| ^t |d\xi |\asymp \int _{C_r} |h'(z)|^t |dz| \quad \text{where} \quad C_r = \exp (\sg )\,.
\eeq
The conclusion comes now from Formula \ref{x.21}, in fact from Proposition \ref{y.2}.
In order to be able to apply it notice that there exists $c>0$ such that $\diam (C_r) \geq c$ for every $1<r<2$ since $C_r= \exp \circ G_T (\ga )$, since all the maps $\exp\circ G_T$, $t\ge 1$, are uniformly bi--Lipschitz
and since $\diam (\ga )\geq 1$. Therefore, 
\beq\label{new 22}
\b_\infty (t)=\b_h (t)= t-1+ \frac{\P(t)}{\log d}\,.
\eeq
The behavior of the pressure function is perfectly understood thanks to  \cite{PRS03} and \cite{PR11}. 
In particular, the hyperbolic dimension $\HypDim(p)$ is the first zero of the pressure function $\P$ (see \cite{PRZ99})
and, having this formula for the hyperbolic dimension, Item (2) then follows from Zdunik's work \cite{Zdunik90}.
\epf

\subsection{Linearizers of TCE polynomials}
It is known that the attracting basin of infinity of a polynomial has nice geometry as long as the polynomial has some expansion.
 Carleson, Jones and Yoccoz \cite{CJY94} have shown that $A_p(\infty )$ is a John domain if and only if $p$ is semi-hyperbolic.
 Graczyk and Smirnov  \cite{GS98} considered Collet-Eckmann rational functions. Their result states that attracting and
 super-attracting components of the Fatou set are H\"older if and only if the function is Collet-Eckmann. There is a useful concept which captures the essential features of Collet-Eckmann maps, the one of Topological Collet-Eckmann rational functions. There are various characterizations of such functions. Several of them have been provided in the paper \cite{PRS03} by Przytycki, Rivera-Letelier and Smirnov.  A partial version of their results is this.

\bthm[\cite{GS98}, \cite{PRS03}] \label{x.16}
Let $p:\oc\to\oc$ be a polynomial and let $A_p(\infty)$ be its attracting basin of infinity. Then, the following conditions are equivalent:

\begin{enumerate}
\item The polynomial $p$ is TCE,

\sp\item $A_p(\infty)$ is a H\"older domain,

\sp\item Negative pressure: $\P(t)<0$ for large values of $t$.
\end{enumerate}
\ethm
 
\fr For more about characterizations of TCE maps see \cite{GS98, PRS03}. It turns out that we get an additional characterization
of the TCE property:

\bthm\label{x.6}
Let $p:\oc\to\oc$ be a polynomial with connected Julia set, let $z_0$ be a repelling fixed point of $p$ and let $f:\C\to\C$ be a corresponding linearizer such that $f\in \cD$ (see Theorem~\ref{x.6 general}). Then, the following are equivalent:
\begin{itemize}
\item[(a)] $p$ is Topological Collett--Eckmann polynomial (equivalent to $A_p(\infty )$ being a H\"older domain).

\item[(b)] All the connected components of $f^{-1}(\D^*)$ are H\"older tracts.

\item[(c)] $f$ has negative spectrum. 
\end{itemize}
\ethm

\bpf
Equivalence of (a) and (c): The pressure function $t\mapsto \P(t)$ of $p$
is always convex and decreasing and has a first zero which is $\HypDim(p)$. Theorem \ref{x.16}, more precisely Section 4 in \cite{PRS03},  tells us
that $p$ is TCE if and only if the pressure function is strictly decreasing and then $\HypDim(p)$ is the only zero.
Hence, given the relation  \eqref{new 22} between the integral mean spectrum of $f$ and the pressure $\P$ of $p$, it follows
that $f $ has negative spectrum if and only if $p$ is TCE. 

Proposition \ref{28} shows that (b) implies (c). It remains to proof that (a) implies (b), i.e. that the maps $\ph\circ T$ satisfy the H\"older condition  \eqref{holder 1'}  with uniform constants.
In order to do so, we must estimate the derivative $|(\ph \circ T)'(1)|$. Equality \eqref{z.12} 
combined with \eqref{new 31}
show that $|(\ph \circ T)'(1)|\asymp |\l |^N |H'(z_T)|$ where $z_T=G_T(1)$. If $L$ is the bi-Lipschitz constant of $G_T$ then
$z_T \in \D(z_T, 1/L )\subset \cH$, hence $\Re z_T \geq 1/L$. On the other hand, by the choice of $n_0$ hence of $N$, we have $z_T\in V''=H^{-1}(V'\cap \cH )$ and on $V"\cap \{\Re z \geq 1/L\}$
we certainly have $|H'|\asymp 1$. Thus 
$$|(\ph \circ T)'(1)|\asymp |\l |^N |H'(z_T)|\asymp |\l |^N.$$
From the expression  \eqref{new 30} we get
\beq\label{x.10}
\ph\circ T =\l^N\circ \ph _{|Q_{1,T}} \circ H\circ G_T \quad \text{on} \quad Q_1\,.
\eeq
We know that $\ph _{|Q_{1,T}}$ and $G_T$ are uniformly bi-Lipschitz. The map $H$ is evaluated on points of
the set $G_T(Q_1)$ which is a subset of $V''$.
But $H_{|V''}$ is a H\"older map since now the polynomial $p$ is TCE and thus the conformal map $h:\D^*\to A_p(\infty)$ is H\"older.
 Its H\"older exponent is $\a$ and denote by $\Ga$ its H\"older constant. 
  Then, making use of  \eqref{x.10}, we get for all points $z_1,z_2\in Q_1$ that
$$
|\ph\circ T(z_1)-\ph\circ T(z_2)|\asymp |\l|^N |g_T(z_1)-g_T(z_2)|\preceq |(\ph\circ T)'(1)| \,\Ga\, |z_1-z_2|^\al
\,.
$$
\epf


\subsection{Poincar\'e functions without negative spectrum} \label{without ns}
 For functions with negative spectrum, the series defining the transfer operator converges exponentially fast. In general, i.e. without assuming negative spectrum, this series can still converge. We illustrate this here by considering arbitrary Poincar\'e functions, i.e. also those without negative spectrum. They are very special entire functions because of the functional equation \eqref{x.3}. This equation allows us to do direct calculations, in the spirit of  \cite{ER15},  in order to relate $\pft \1 (w)$ to the Poincar\'e series 
 $$
\cP_t(\xi) := \sum _{N\geq 1} \sum_{\eta \in p^{-N}(\xi )}\big|\left( p^{N} \right)'(\eta )\big| ^{-t} 
$$
 of the polynomial $p$ evaluated at the point $\xi$. Let $\d_{cr}$ be the critical exponent of this series  so that 
 $\cP_t(\xi) <\infty $ for $t>\d_{cr}$ and  $\cP_t(\xi) =\infty $ for $t < \d_{cr}$, then we have the following.

\bthm\label{irr pfunctions} 
Let $f\in \cD$ be a  linearizer of a polynomial $p:\oc\to\oc$ with connected Julia set such that $\cJ_f \subset A_p(\infty )$.
Then, there exists a neighborhood $\cV$ of the Julia set $\J_f$ such that for every $t>\tstar=\d_{cr}$  the Perron--Frobenius operator $\pft$ is well defined and bounded on $\cC_b (\cV )$. Moreover,
$$ 
\pft \1 (w) \asymp ( \log |w|)^{1-t} \quad \text{for every} \quad w\in \cV,
$$
with comparability constant depending on the whole initial data.
\ethm

\bpf 
Let $h:\D^*\to A_p(\infty )$ be the Riemann map such that
\beq\label{z.2}
h(z^d)=p\circ h(z) \; \text{ in }\; z\in\D^*\,.
\eeq
Since $\cJ_f \subset A(\infty )$, there exists $R_0>1$ such that $\J_f\subset h(\D_{R_0}^*)$. Denote 
$$
G_0 := h\(A(R_0^{1/d},R_0)\)\,.
$$
This is a fundamental annulus for the action of $p$ in $A_p(\infty)$. We also need such an annulus for the action of 
$p$ near the repelling fixed point $z_0$:
$$ 
V_0 := f\(A(r,|\l|r)\),
$$
where $r>0$ is so small that $f$ is univalent on $\D_{2| \l | r}$ and 
\beq \label{z.1}
f\(\D_{2| \l | r}\)\subset \C \setminus h\(\D_{R_0^{1/d}}^*\)\,.
\eeq
Notice that $V_0\cap \J_p \neq \emptyset$ since otherwise $z_0$ would be an isolated point of $\J_p$. Therefore there exists 
$M\geq 1$ such that 
\beq\label{z.7}
p^M(V_0)\supset h\(A(1,R_0)\)\,.
\eeq
Let in the following $w$ be an arbitrary point of 
$$ 
\cV:=h\(\D_{R_0^{1/d}}^*\)\supset \J_f.
$$ 
Then there exists a unique integer $N_w\ge 0 $ such that $w\in p^{N_w}(G_0)$. It then follows from \eqref{z.2}, iterated $N_w$ times, that 
 \beq\label{z.3} 
 d^{N_w} \asymp \log |w| \,.
 \eeq
In order to estimate $\pft\1 (w)$ we have to estimate $|f'(z)|_1$ for all $z\in f^{-1}(w)$. If $z$ is such a pre-image, i.e. if $z\in f^{-1}(w)$, then there exists a unique integer $n\geq 1$ such that 
$$
\l^{-n}z\in A(r,|\l|r).
$$ 
Then $w= f(z)= p^n\circ f \circ \l^{-n}(z)$. By \eqref{z.1}, $n> N_w$. Setting $N:=n-N_w$, $\eta := f(\l^{-n}z)\in V_0$ and $\xi :=p^N(\eta )\in G_0$, we get
 \beq \label{z.4}
 f'(z)= \left( p^{N_w} \right)'(\xi ) \left(p^N \right)'(\eta ) f'(\l^{-n} z) \l^{-n}\,.
 \eeq
Since $f$ is univalent on $\D_{2| \l | r}$ and since 
 $\lambda^{-n}z\in A(r,|\l|r)\sbt\D_{|\l | r}$, we have $|f'(\l^{-n} z)|\asymp 1$.
The factor $\left( p^{N_w} \right)'(\xi )$ can be estimated as follows. Since $|h'|\asymp 1$ on $\ov{\D_{R_0^{1/d}}^*}$, we have $|h^{-1} (w)|\asymp |w|$,  $|(h^{-1})'(\xi )| \asymp 1$ and $|h'\big( a^{d^{N_w}}\big)|\asymp 1$
where $a=h^{-1}(\xi)$. Therefore,
$$
\big|\left( p^{N_w} \right)'(\xi ) \big|= \big|h'\big( a^{d^{N_w}}\big)\big| d^{N_w} |a|^{d^{N_w}-1} \big|(h^{-1})'(\xi ) \big|
\asymp  \frac{d^{N_w}}{|a|} \big| h^{-1} (w)\big| \asymp d^{N_w} |w|\,.
$$
Inserting this into \eqref{z.4} leads to
$$
|f'(z)|_1= \Big| \frac{f'(z)}{w}z\Big|\asymp d^{N_w} \big|\left( p^{N} \right)'(\eta )\big| |\l^{-n}z| \asymp  d^{N_w}\big| \left( p^{N} \right)'(\eta )\big|.
$$
Finally this gives
\beq\label{z.6}
\pft\1 (w) \asymp \sum_{\substack{\xi\in G_0 \\ p^{N_w}(\xi)=w}} \sum _{N\geq 1} \sum_{\eta \in p^{-N}(\xi )\cap V_0} \left( d^{N_w} \big|\left( p^{N} \right)'(\eta )\big|  \right)^{-t} \,.\eeq
Let 
$$
\cP_t(\xi , V_0) := \sum _{N\geq 1} \sum_{\eta \in p^{-N}(\xi )\cap V_0}\big|\left( p^{N} \right)'(\eta )\big| ^{-t}
$$ 
and remember that $
\cP_t(\xi) $ designs the corresponding full Poincar\'e series of the polynomial $p$ evaluated at the point $\xi$. 

\bclaim \label{z.5} There exists a constant $c_t >0$ such that
$$
c_t \cP_t(\xi )\leq \cP_t(\xi , V_0) \leq \cP_t(\xi )
$$
for all $\xi \in G_0$.
\eclaim

\bpf Recall that the integer $M$ has been introduced in \eqref{z.7}.
We have
$$
 \cP_t(\xi , V_0)  \geq \sum_{N> M} \sum_{\eta \in p^{-N}(\xi )\cap V_0}\big|\left( p^{N} \right)'(\eta )\big| ^{-t}.
$$
Note that for every integer $k>0$ we have that $p^{-k} (G_0)\subset A(1,R_0)$ and thus, using \eqref{z.7}, we conclude that
for every $z\in p^{-k}(\xi )$, $\xi\in G_0$, there exists at least one point $\eta \in p^{-M} (z)\cap V_0$.
Therefore,
$$
 \cP_t(\xi , V_0)  
 \geq \inf_{\eta \in V_0} \big|\big( p^{M} \big)'(\eta )\big| ^{-t}
 \sum_{k> 0}  \sum_{z \in p^{-k}(\xi )\cap V_0}\big|\big( p^{k} \big)'(z )\big| ^{-t}
 =c_t \cP_t(\xi),
 $$
 where $c_t:=\inf_{\eta \in V_0} \big|\big( p^{M} \big)'(\eta )\big| ^{-t}$.
The other inequality in Claim \ref{z.5} trivially holds and so its proof is complete.
\epf

\fr

Fix arbitrarily $\xi_0\in G_0$. Koebe's Distortion Theorem implies that $\cP_t(\xi )\asymp \cP_t(\xi _0)$ and
thus, by Claim \ref{z.5}, 
\beq\label{1_2017_08_03}
\cP_t(\xi , V_0)\asymp \cP_t(\xi _0)
\eeq for every $\xi \in G_0$. On the other hand, $w$ has exactly $d^{N_w}$
preimages $z\in p^{-N_w}(w)$ and they are all in $G_0$. We can therefore deduce from \eqref{z.6} that
$$
\pft\1 (w) 
\asymp  d^{-t N_w}\sum_{\substack{\xi\in G_0 \\ p^{N_w}(\xi)=w}}\cP_t(\xi _0, V_0)  \asymp d^{N_w (1-t)} \cP_t(\xi _0)
\asymp d^{N_w (1-t)} \,.
$$
The conclusion now follows directly by applying  \eqref{z.3}.
\epf


\smallskip

\section{The Classics of Thermodynamic Formalism: \\ Conformal Measures and Beyond}\label{section MD}

Let $f\in \cD$ be a function with negative spectrum. Then the whole 
thermodynamic formalism can be established for $f$, word by word,
exactly as it was done in \cite{MyUrb08, MUmemoirs} except for \cite[Lemma 5.13]{MUmemoirs} which is the key point in the
construction of conformal measures. Since we provide below, in Proposition \ref{x.1}, a proof of this missing point, 
we finally show that all the relevant results comprising Thermodynamical Formalism, the ones established in \cite{MUmemoirs} 
 and stated below, hold. Combined with Theorem \ref{56} this shows 
 Theorem \ref{thm intro}.

In Section \ref{without ns} we considered some special class, of entire functions, called Poincar\'e functions, that do not necessarily have negative spectrum. As Theorem \ref{irr pfunctions} shows, these perfectly satisfy the assumptions of Proposition \ref{x.1}. 
Consequently, this proposition and all the results, Theorem~\ref{theo main X} through Theorem~\ref{thVariational Principle}, of the present section are also valid for these functions.

\

\sp  
\fr
- The Perron-Frobenius-Ruelle Theorem \cite[Theorem 5.15]{MUmemoirs}.

\bthm\label{theo main X}
If $f\in \cD$ is a function with negative spectrum and $t>\tstar$, then the following are true.

\sp\begin{itemize}
    \item[(1)] The topological pressure $\P(t )=\lim_{n\to\infty}
    \frac{1}{n} \log \pf^n_t\1(w)$ exists and is independent of $w\in \jul_f$.
    
\item[(2)] The function $(\tstar,+\infty)\ni t\longmapsto\P(t)\in\R$ is convex, thus continuous, in fact real--analytic, strictly decreasing, and $\lim_{t\to+\infty}\P(t)=-\infty$.
    
    \item[(3)] There exists a unique $\l|f'|_1^t$--conformal
    measure $m_t$ and necessarily $\l=e^{\P(t)}$. Also, there exists a unique Gibbs state $\mu_t $, i.e.
    $\mu_t$ is $f$-invariant and equivalent to $m_t$.
    \item[(4)] Both measures $m_t$ and $\mu_t$ are ergodic and supported on the radial (or conical) Julia
    set $\rad(f)$.
    \item[(4)] The density $\den_t :=d\mu_t/dm_t$ is an everywhere positive
    continuous and bounded function on the Julia set $\jul_f$.
\end{itemize}
\ethm

\

\fr
- The Spectral Gap \cite[Theorem 6.5]{MUmemoirs}

\bthm\lab{t6120101}
If $f\in \cD$ is a function with negative spectrum and $t>\tstar$, then the following are true.

\sp\begin{itemize}
\item[(a)] The number $1$ is a simple isolated eigenvalue of
the operator $\hat\pf_t:=e^{-\P(t)}\pf_t:\H_\b\to \H_\b$ ($\b\in(0,1]$ is arbitrary and $\H_\b$ is the Banach space of real--valued bounded H\"older continuous defined on $\jul_f$) and all other eigenvalues are contained
in a disk of radius strictly smaller than $1$.
\item[(b)] There exists a bounded linear operator $S:\H_\b\to \H_\b$ such that
$$
\hat\pf_t=Q_1+S,
$$
where $Q_1:\H_\b \to \C \den$ is a projector on the eigenspace $\C
\den$, given by the formula 
$$
Q_1(g)=\lt(\int g \, dm_\phi\rt)\den_t,
$$
$Q_1\circ S=S\circ Q_1=0$ and
$$
||S^n||_\b\le C\xi^n
$$
for some constant $C>0$, some constant $\xi\in (0,1)$ and all $n\ge 1$.
\end{itemize}
\ethm

\

\fr
- \cite[Corollary 6.6]{MUmemoirs}

\bcor \label{5.1.2}
With the setting and notation of Theorem \ref{t6120101} we have, for every $n\geq 1$, that $\npf ^n = Q_1 +S^n$
and that $\npf ^n (g)$ converges to $\left( \int g\, dm_\phi \right) \den $ exponentially fast when $n\to \infty$. Precisely,
$$\lt\| \npf ^n (g) -\left( \int g\, dm_\phi \right) \den \rt\| _\b =\| S^n (g) \| _\b \leq C \xi ^n \|g\| _\b \quad , \;\; g\in H_\b.$$
\ecor

\

\fr
- Exponential Decay of Correlations \cite[Theorem 6.16]{MUmemoirs} 

\bthm\lab{t3120301}
With the setting and notation of Theorem \ref{t6120101} there exists a large class of functions $\psi_1$ such that for all $\psi_2\in L^1(m_t)$ and all integers $n\ge 1$, we have that
$$
\lt|\int(\psi_1\circ f^n\cdot\psi_2)\,d\mu_t- \int\psi_1\,d\mu_t\int\psi_2\,d\mu_t\rt|\le
O(\xi^n),
$$
where $\xi\in(0,1)$ comes from Theorem~\ref{t6120101}(b), while the big ``O'' constant depends on both $\psi_1$ and $\psi_2$.
\ethm

\

\fr
- Central Limit Theorem \cite[Theorem 6.17]{MUmemoirs}

\bthm
With the setting and notation of Theorem \ref{t6120101} there exists a large class of functions $\psi$ such that the sequence of random variables
$$
\frac{\sum_{j=0}^{n-1}\psi\circ f^j-n \int\psi\,d\mu_t}{\sqrt{n}}
$$
converges in distribution with respect to the measure $\mu_t$ to the Gauss (normal) distribution ${\mathcal N}(0,\sigma^2)$ with some $\sg>0$. More precisely,
for every $t\in \R$,
$$
\begin{aligned}
\lim_{n\to\infty}\mu_t\bigg(\bigg\{z\in \jul_f: \frac{\sum_{j=0}^{n-1}\psi\circ f^j(z)-n \int\psi\,d\mu_t}{\sqrt{n}}&\le t\bigg\}\bigg)= \\
&= {1\over \sigma \sqrt{2\pi}}\int_{-\infty}^t\exp\lt(-\frac{u^2}{2\sigma^2}\rt)\, du.
\end{aligned}
$$
\ethm

\

\fr
- Variational Principle \cite[Theorem 6.25]{MUmemoirs}
\bthm \label{thVariational Principle}
If $f\in \cD$ is a function with negative spectrum and $t>\tstar$, then the
$f$--invariant measure $\mu_t$ is the only equilibrium state of the
potential $-t\log|f'|_1$, that is
$$
\P(t)=\sup\lt\{{\rm h}_\mu(f)-t\int_{\jul_f}\log|f'|_1\, d\mu\rt\},
$$
where the supremum is taken over all Borel probability $f$-invariant
ergodic measures $\mu$ with $\int_{\jul_f}\log|f'|_1\, d\mu>-\infty$, and
$$
\P(t)={\rm h}_{\mu_t}(f)-t\int_{\jul_f}\log|f'|_1\, d\mu_t.
$$
\ethm

We will obtain conformal measures following the approach in \cite[Section 5.3]{MUmemoirs}. For dynamical systems with compact Julia set these measures can be produced either as fixed points of dual Perron-Frobenius operators or as weak limits of some atomic measures using the fact
that the space of probability measures on the Julia set is weakly compact. In the present setting the Julia set $\J_f$ is an unbounded subset of $\C$
and so the key point is to establish the tightness (in $\C$) of an appropriate sequence of measures. This can be done by following  \cite[Section 5.3]{MUmemoirs} since we have the following analogue of  \cite[Lemma 5.13]{MUmemoirs}:

\bprop\label{x.1}
Let $f\in \cD$ and suppose that there exists $t>0$ for which $\pft$ is a bounded operator of $\cC_b(\ov \Om )$ with
$$
\lim_{w\to \infty \; , \;  w\in \Om} \pft \1 (w) =0\;.
$$
 Then 
$$
\displaystyle \lim_{S\to \infty} \|\pft  \1_{\D_S^*}\|_\infty=0\,.
$$
\eprop

\bpf
Let $\ep >0$ and let  $R_\ep >e^\ga$ be so large that 
$\pft \1 (w) < \ep$ for all $w\in \Om$ with $|w| > R_\ep$.
Then clearly $\pft  \1_{\D_S^*} (w) \leq \ep $
for every $S>0$ and every $w\in\Om$ with $ |w| >R_\ep$.

We are left to consider points 
$$
w\in\ov\Om \  \ {\rm  with } \   |w| \leq R_\ep.
$$
The set $\ov \Om\cap \ov\D_{R_\ep}$ is compact and thus,
for every $\d >0$, it admits a finite covering by $\d$--disks with centers $w_1,...,w_N$ in 
$\ov \Om$. Because of the disjoint assumption (see  \eqref{1}) we can 
choose $\d>0$ such that
every disk centered in $\ov \Om$
and of radius $2\d$ does not intersect the singular set $S(f)$. 
Then all inverse branches of all iterates of $f$ can be defined on the disks
of the $\d$--covering and they have bounded distortion.

Consider an arbitrary point 
$
w\in\ov\Om 
$ and let $j\in \{1,...,N\}$ such that $w\in \D(w_j ,\d )$. Since all inverse branches of $f$ are
well defined on $\D(w_j ,\d )$ we have natural pairings of inverse images of $w$ and $w_j$.
Indeed, if $z\in f^{-1}(w)$ then there exists an inverse branch $f_z^{-1}$ of $f$ defined on the
disk $\D(w_j ,\d )$ such that $z= f_z^{-1} (w)$. To this we naturally associate $z_j\in f^{-1}(w_j)$
defined by $z_j= f_z^{-1}(w_j)$ and, the disjoint
type assumption implying expansion \cite{Rempe09},  we have
$$
|z-z_j|\leq c|w-w_j|\leq c\d
$$
for some constant $c>0$ depending on $f$ only. In particular,  if $|z|> S$ then $|z_j|> S-c\d$.

Let from now on now $S> c\d$ so that $S'= S-c\d >0$.
The previous considerations along with the bounded distortion property imply that there exists $K<\infty$
such that
$$
\pft  \1_{\D_S^*} (w) \leq K\pft  \1_{\D_{S'}^*} (w_j) \quad \text{for every} \quad w\in \D(w_j, \d ) \;, \; j=1,...,N\,.
$$
By its very definition,  $\pft  \1_{\D_S^*} (w)$ takes into account only the preimages $z\in f^{-1} (w)$ for which $|z|\geq S$. 
Since  $\pft \1 (w) $ is convergent, for every $w\in \Om$, there exists $S>c\d$ such that 
$$
\pft  \1_{\D_{S'}^*} (w_j)<\ep/K\quad \text{for every $j=1,...,N$.}
$$
 Hence,
$
\pft \1_{\D_S^*} (w)<\ep
$
for all $w\in \ov \Om\cap \ov\D_{R_\ep}$.
This completes the proof of Proposition \ref{x.1}.
\epf

\smallskip

\section{Thermodynamics: Bowen's Formula}\label{section MD'}

Let $f\in \cD$ have negative spectrum.
The pressure function
introduced in the previous section along with its properties established in Theorem~\ref{theo main X} (2) allows us to provide a closed formula for the Hausdorff dimension of the radial Julia set of $f$. This quantity is called the \emph{hyperbolic dimension} of $f$, is denoted  by $\HypDim(f)$ and is also known (see \cite{Rempe09-hyp}) to be the supremum of the Hausdorff dimension of all hyperbolic subsets of $\cJ_f$.
Here is a reformulation of Theorem \ref{Bowen's}.

\bthm[Bowen's Formula] \label{42}
Let $f\in \cD$ have negative spectrum. Then, the function $(\tstar,+\infty)\ni t\mapsto \P(t)$ has a (unique) zero 
$h>\tstar$
if and only if $\HypDim(f)>\tstar$.
In this case we have 
$$
\HypDim(f) = h\,.
$$
\ethm

\bpf
Since $f\in \cD$ is a function with negative spectrum the thermodynamic formalism of Section~\ref{section MD} applies. In particular, for every $t>\tstar$ there exists an $e^{\P(t)}|f'|_1$--conformal measure. If for some $h>\tstar$ we have $\P(h)=0$, then the corresponding conformal measure is frequently called geometric conformal measure, i.e. $e^{\P(h)}=1$. The proof of Theorem~1.2 in \cite{MyUrb08} then applies yielding $\HypDim(f) = h(>\tstar)$.

Conversely, if $\HypDim(f)>\tstar$, then (see \cite{Rempe09-hyp}) there exists a hyperbolic (compact) set $X\sbt \cJ_f$ such that $\HD(X)>\tstar$. Then, see  \cite{PUbook}  for example, $\P(f|_X,-\HD(X)\log|f|_X'|)=0$. Therefore, $\P(\HD(X))\ge 
\P(f|_X,-\HD(X)\log|f|_X'|)=0$. In conjunction with Theorem~\ref{theo main X} (2) this implies that there exists $h\ge \HD(X)$ such that $\P(h)=0$. As $\HD(X)>\tstar$, the proof is complete.
\epf

The issue raised by Theorem~\ref{42} is to be able to tell whether $\HypDim(f)>\tstar$. We know that there is a quite general class of functions for which this holds. Indeed, Bara{\'n}ski, Karpi{\'n}ska and Zdunik showed in \cite{BKZ-2009} that $\HypDim(f)>1$ for every function $f\in \cD$. Along with Theorem~\ref{42} this implies the following.

\bprop
Let $f\in \cD$ have negative spectrum with  $\tstar\leq1$. Then, the pressure function  $(\tstar,+\infty)\ni t\mapsto \P(t)$  of $f$ has a unique zero, call it $h$, and 
$$
\HypDim(f) = h\,.
$$
\eprop

The problem when $\tstar >1$ is that then the $\tstar$--Hausdorff measure of the boundary of a tract may be zero.  If this is quantitatively not the case in the sense that $\ph_T (I)$ has $\tstar$--measure greater than some strictly positive constant
and if the tract is H\"older
then the hyperbolic dimension can be estimated like in \cite{BKZ-2009}. This has been observed and worked out for a family of examples in Proposition 4.3 of \cite{My-HypDim}.
The model functions of the Proposition 4.3 in \cite{My-HypDim} have all the property $\HypDim (f) > \tstar$.

\medskip 

In general we can use our optimal estimates for the transfer operator in order to 
show directly that the pressure function  has a zero. Here, it is done for a large class of linearizers. Let us mention again that \cite{DR} contains a more general version of this result.

\bprop\label{38}
Let $p:\oc\to\oc$ be a hyperbolic polynomial with connected Julia set. Let $z_0\in \J_p$ be a repelling fixed point of $p$
and let $f$ be an associated linearizer of disjoint type. Then, the pressure  function $(\tstar,+\infty)\ni t\mapsto \P(t)$ of $f$ has a zero which we denote by $h$. In consequence,
$$
\HypDim(f) = h\,.
$$
\eprop

\bpf 
We are to estimate $\pft \1$ for $t>\tstar$ near $\tstar$. 
It suffices to show that there exists $R>1$ such that 
\beq \label{43}
\pft \1_{\D_R} (w) \geq 2 \quad \text{for all}\quad w\in \D_R\cap \J_f
\eeq
since then it follows by induction that $\pft ^n\1 (w) \geq 2^n$ for all $w\in  \D_R\cap \J_f$; thus that $\P(t)\ge \log 2>0$. This, along with Theorem~\ref{theo main X} (2) entails the existence of a unique zero $h>\tstar$.

For the special type of functions $f$ we consider here we have the estimate from Theorem~\ref{irr pfunctions}:
$$
\pft \1 (w) \asymp  (\log |w|)^{1-t}.
$$
With the notations of the proof of Theorem \ref{irr pfunctions}, let
$$
A_n:= \{z\in \C: |\l|^n r < |z| \leq |\l|^{n+1} r  \}\;\; \text{ and } \;\; R_n=|\l|^{n +1}r \; , \;\;  \quad n\geq 0\,.
$$
Then for every integer $M\ge 1$ we have that
$$
\pft (\1_{\D_{R_M}}) (w)= \sum_{n = 0}^M \sum_{z\in f^{-1}(w)\cap A_n} |f'(z)|_1^{-t}\; .
$$
Because of \eqref{z.6} and  \eqref{z.3} we see that the sum over the preimages of $w$ lying in $A_n$ is approximately 
\beq\label{1_2017_08_12}
(\log |w|)^{1-t}\sum_{\eta \in p^{-N}(\xi )\cap V_0}\big|\left( p^{N} \right)'(\eta )\big| ^{-t}
\eeq
where $N=n-N_w$. Since the polynomial $p:\oc\to\oc$ is hyperbolic, it is well known, in fact this is a theorem of thermodynamic formalism of distance expanding dynamical systems (see \cite{PUbook} for ex.) that the sum of \eqref{1_2017_08_12}
is approximately $\l_{t,p}^N$ where $\log\l_{t,p}$ is the topological pressure of $p$ 
evaluated at $t$. Also, $t>\tstar=\HypDim(p)=\HD (\J_p)$ because of Theorem \ref{x.6} and 
again since $p$ is hyperbolic. In particular, $\l_{t,p}<1$ for all $t> \tstar$ and $\lim_{t\to \tstar}\l_{t_,p}=\l_{\tstar,p} =1$.
It follows that
$$
\pft (\1_{\D_{R_M}}) (w)\asymp (\log |w|)^{1-t} \frac{1-\l_{t,p}^{M+1}}{1-\l_{t,p}}\asymp  (\log |w|)^{1-t}  M,
$$
where the latter comparability holds if we make the particular choice $ M=M(t):=\frac2{\log \l_{t,p}}$. So, if $w\in \D_{R_M}\cap \J_f$ then
$$
\pft (\1_{\D_{R_M}}) (w) \succeq M(t)^{2-t}.
$$
Since $\lim_{t\downto \tstar}M(t)=+\infty$, we immediately see that there exists $t>\tstar$ such that \eqref{43} holds with $R=R_{M(t)}$.
The proof is complete.
\epf

\bpf[Proof of Theorem \ref{analyticity}] This theorem is an immediate consequence of 
Theorem~\ref{x.6 general}, Corollary \ref{38}, Corollary~\ref{5.1.2}, Theorem~\ref{42} and Theorem~9.11 of \cite{MUZ} (our present theorem is a very special case of Theorem~9.11 as we now consider only deterministic systems).
\epf

\section{Functions of Class $\cS$} \label{section S}

We finally consider functions of the Speiser class $\cS$. For this more restrictive class
the present theory of 
thermodynamic formalism
can be extended in a straightforward way
 to hyperbolic functions. This section contains also the promised second application of quasiconformal invariance
of H\"older tracts in Proposition  \ref{w.3} and provide a prove for Theorem \ref{thm 1.4}
of the Introduction.

\subsection{Hyperbolic Functions in Class $\cS$}\label{dishyp}
The object here is to explain how to pass from the disjoint type case to hyperbolic functions. In order to do so,
let us consider a function $f$ having negative spectrum and the properties of class $\cD$ except for the disjoint type property. Instead, $f$ is assumed to be hyperbolic and of class $\cS$.

We assume as usual that $S(f)\subset \D$, that $0\in \F(f)$ and we fix arbitrarily  $\ga >0$.
Then, exactly as in the disjoint type case, the conclusions of Theorem~\ref{50} hold for every $|w| > e^\ga $. 
We are thus to consider only points 
$$
w\in \ov \D _{e^\ga}\cap \J_f.
$$ 
We will compare $\pft \1 (w)$ to $\pft \1 (\xi)$
where $\xi\in \C\setminus \ov \D_{e^\ga}$ is an arbitrarily fixed point. This goes exactly as in \cite[Section 10]{Bishop-S-2016} (and this is the only point where class $\cS$ rather than merely class $\cB$ is needed) by employing the bounded distortion argument.  Indeed, it suffices to connect $w$ to $\xi$ by a piecewise smooth path
$\sg$ of Euclidean length uniformly bounded above with respect two points  $w\in \ov \D _{e^\ga}$ and such that for some fixed $\d>0$, 
the $\d$--neighborhood of $\sg$ does not intersect $S(f)$. Let us recall from \cite[Section 4.2]{MUmemoirs} that there exists good distortion estimates
for $|(f^n)'|_1$. In conclusion, all of this  shows that Theorem \ref{thm intro} holds for $f$.


\subsection{Analytic families of class $\cS$}
Two entire functions $f$ and $g$ are (topologically) equivalent if there exist two homeomorphisms $\Phi , \Psi:\C\to\C$
such that 

\beq\label{w.1}
\Psi \circ g = f\circ \Phi \;.
\eeq
Given $g\in \cS$, Eremenko and Lyubich \cite{EL92} showed that the set $\cM_g$ of all functions $f\in \cS$ equivalent to $g$
has a natural structure of a complex analytic manifold. It can be parametrized with the help of the singular values $\{a_1,...,a_q\}=S(g)$ of $g$.

\bprop\label{w.3}
Let $g\in \cS$ be a function having finitely many tracts all of which are H\"older. Then all tracts of every function $f\in \cM_g$ are H\"older, and thus
all functions of $\cM_g$ have negative spectrum.
\eprop

\fr The proof of this fact relies on a special choice of homeomorphisms in the equivalence relation \eqref{w.1}.
In fact, they can be freely chosen in an isotopy class without changing $f,g$ and, in particular, there is a quasiconformal choice for these homeomorphisms (see again \cite{EL92}). We first shall prove the following.

\blem\label{w.2}
Let $g\in \cS$. If $f\in \cM_g$, then the homeomorphism $\Psi$ in \eqref{w.1} can be chosen quasiconformal and such that, 
for some $R\geq 1$,
$$
\Psi|_{\D_R^*}\equiv \Id|_{\D_R^*}\,.
$$
\elem

\bpf
Assume without loss of generality that $S(g)=\{a_1,...a_q\}\subset \D$ and let $\psi: \C\to \C$ be a quasiconformal homeomorphism such that  \eqref{w.1} holds for the given maps $g\in \cS$ and $f\in \cM_g$.
A standard application of the Ahlfors-Bers-Bojarski Measurable Riemann Mapping Theorem is that $\psi^{-1}$ can be embedded into 
a holomorphic motion of quasiconformal mappings $\psi_\l^{-1}:\C\to\C$, $\l\in \D$. In particular, $\psi_0^{-1}=\Id$ and 
$\psi_{\l_0}^{-1}=\psi^{-1}$ for some $\l_0\in \D$. Define $R$ to be any number
$$
\ge 2 \max\big\{\diam (\psi _\l (\D)):|\l|\leq |\l_0|\big\},
$$
and consider the holomorphic motion
$
z\mapsto z_\l
$
defined on $\psi(\D \cup \D^*_R)$ by 
$$
\text{$z_\l = z$ if $z\in \psi (\D)$ \  and $\  z_\l= \psi_\l ^{-1} (z)$ if $z\in \psi (\D_R^*)$.}
$$
By Slodkowski's version of Ma\~n\'e-Sad-Sullivan's $\l$-Lemma, this holomorphic motion has 
an extension to a holomorphic motion $h_\l :\C\to\C$, $|\l|\leq |\l_0|$. Then $h_\l$, thus also
$$
\Psi_\l := h_\l \circ \psi \quad , \quad \l \in \D \,,
$$
 are quasiconformal mappings. The map we look for is $\Psi_{\l_0}$ and $t\mapsto \Psi_{t\l_0}$, $t\in [0,1]$ is an isotopy between $\psi = \Psi_0$
and $\Psi_{\l_0}$ that does not move the points $a_j$, $j=1,...,q$; in fact it does not move any point of $\D$. 
It suffices now to apply Proposition 2.3 in \cite{ER15} since it shows that there exists $\Phi_{\l_0}:\C\to\C$ quasiconformal such that  \eqref{w.1} holds with $\Psi$ and $\Phi$ replaced respectively by $\Psi_{\l_0}$ and $\Phi_{\l_0}$.
\epf

\bpf[Proof of Proposition \ref{w.3}]
Let $g,f$ be as in  Proposition~\ref{w.3} and let $R>1$ so large that the assertion of Lemma~\ref{w.2} holds and that all the components of $g^{-1}(\ov \D_R^c)$ are H\"older.
Then, if $\Psi$ is given by Lemma \ref{w.2} and if $\Phi$ is the corresponding quasiconformal map such that \eqref{w.1} holds, then clearly $\Phi$ identifies the components of 
$g^{-1}(\ov \D_R^*)$ with those of $f^{-1}(\ov \D_R*)$. Proposition~\ref{w.3} follows now from Lemma~\ref{w.4}.
\epf

\subsection{Proof of Theorem \ref{thm 1.4} and of Corollary  \ref{thm 1.4'}} 
If $g$ is  as in Theorem \ref{thm 1.4}, then Proposition \ref{w.3} yields that every
$f\in \cM_g$ has H\"older tracts and negative spectrum. Therefore, Theorem \ref{thm intro} applies first to any disjoint type map of $\cM_g$ and 
then also to every hyperbolic function of this family $\cM_g$ because of the argument of Section \ref{dishyp}.
This proves Theorem  \ref{thm 1.4}. \hfill $\square$

\medskip

Let now $g\in \cS$ be a linearizer of a polynomial $p$ with connected Julia set. The singular set $S(g)$ equals the post-singular set of the
polynomial $p$ (see \cite{MP12}). Since $g\in \cS$, $p$ must be post-critically finite hence TCE. On the other hand, we may assume
that $g$ is of disjoint type so that it satisfies the conclusion of Lemma \ref{x.36}.
Otherwise it suffices to replace it by $g\circ \ka\in \cM_g$ with sufficiently small $\ka \neq 0$. 
It thus follows from Lemma \ref{x.36} along with Theorem \ref{x.6}
that $g\in \cD$ and that $g$ has finitely many H\"older tracts. It suffices now to apply Theorem \ref{thm 1.4} in order to complete
the proof of  Corollary  \ref{thm 1.4'}. \hfill $\square$


\smallskip
\appendix
\section{Integral means spectrum and pressure} \label{y.1}

For the sake of completeness we  provide here the details related to formula \eqref{x.21} in the setting of the proof of Theorem \ref{x.6}.
Let again $p:\oc\to\oc$ be a polynomial with connected Julia set and let
 $h:\D^*\to A_p(\infty)$ be a Riemann map such that
 \beq\label{z.11}
 h\circ D=p\circ h \quad \text{on} \quad \D^* \quad \text{where}\quad  D(z)=z^d\,.
 \eeq
Suppose we are given a constant $c>0$ and circular arcs $C_r\subset \{|z|=r\}$ with 
 $$
 \diam(C_r)\geq c\quad , \quad r>1\,.
 $$
 Define
  $$\hat \b_{h}(t)=\limsup_{r\to1^+}\frac{\log \int _{C_r} |h'(z)|^t  |dz|}{|\log(r-1)|}\,.$$
Consider also the tree pressure
$$
\P(t,w):=\limsup_{n\to\infty }\frac{1}{n} \log \sum_{z\in p^{-n}(w)}|(p^n)'(z)|^{-t} \quad , \;
t>0 \;\text{ and }\; w\in A_p(\infty )\,.
$$
It has been shown in \cite{PRZ99} that this expression does not depend on $w$. More precisely,
Przytycki has shown that $\P(t,w)$ is the same value for every typical $w\in \C$. Since now the polynomial $p$
is assumed to have connected Julia set, every point of $A_p(\infty )$ is typical in the sense of \cite{PRZ99}.
We can therefore write $\P(t)$ for $\P(t,w)$, for any $w\in A_p(\infty )$.

\bprop\label{y.2} 
If $p:\oc\to\oc$ is a polynomial with connected Julia set
and if $d\geq 2$ denotes its degree, then
$$   
\hat \b_h (t)=      \b_{h}(t)= \frac{\P(t)}{\log d}+t-1\, 
$$
for every $t\ge 0$.
\eprop

We adapt the proof given in \cite{PUbook}. There, the second equality is shown for expanding polynomials.

\bpf
Given $1<r<2$, there exists a unique integer $n\ge 1$ such that $R_r=D^n(r)\in [2,2^d[$. Obviously $r-1\asymp d^{-n}$ and the (arcwise) distance between two consecutive elements of $D^{-n}(R_r)$ is $2\pi r/d^n\comp d^{-n}$. Therefore, applying Koebe's Distortion Theorem, we get the following:
$$
 \int _{C_r} |h'(\xi)|^t  |d\xi | 
 \asymp d^{-n}\sum_{\xi\in C_r\cap D^{-n}(R_r)}  |h'(\xi)|^t.
$$
Iterating the functional equation \eqref{z.11} and taking derivatives gives
$$
h'(R_r) (D^n)'(\xi ) = (p^n)' (h(\xi )) h'(\xi )
$$
for all $\xi\in D^{-n}(R_r)$. Hence, we get for such $\xi$ that  
$$
|h'(\xi)|\asymp d^n |(p^n)'(h(\xi ))|^{-1}\,.
$$
Thus 
\beq\label{41}
  \begin{aligned}
 \int _{C_r} |h'(\xi)|^t  |d\xi | &\asymp d^{n(t-1)}\sum_{\xi\in C_r\cap D^{-n}(R_r)} |(p^n)'(h(\xi ))|^{-t}\\
 &=
d^{n(t-1)} \sum_{z\in h(C_r)\cap p^{-n}(w_r)} |(p^n)'(z)|^{-t}\; , \;\; w_r=h(R_r)\,.
  \end{aligned}
 \eeq
Since the constant $c>0$ from the definition of the circular arcs $C_r$ does not depend on $r>1$, there exists an integer $M\ge 1$ (in fact every integer $M$ large enough is good) such that 
$$
D^M(C_r)= \{|\xi|= r^{d^M}\}
$$ 
for every $r>1$. If $r$ is sufficiently close to $1$ then $n>M$. Then 
$$
h(C_r)\cap p^{-n}(w_r)=h(C_r)\cap p^{-M}(p^{-(n-M)}(w_r)).
$$  
Obviously there exists $R\in(0,+\infty)$ such that 
$$
p^{-(n-M)}(w_r)\sbt \D_R
$$
for all $1<r<2$. Since 
$$
K:=\sup\big\{|(p^M)(v)|:v\in\D_R\big\} <\infty,
$$
we thus have that
$$ \begin{aligned}
\sum_{z\in h(C_r)\cap p^{-n}(w_r)}|(p^n)'(z)|^{-t} &\geq K^{-t}\sum_{\xi\in p^{-{(n-M)}}(w_r)}|(p^{n-M})'(\xi)|^{-t}\\
&\asymp \sum_{\xi\in p^{-{(n-M)}}(h(2))}|(p^{n-M})'(\xi)|^{-t}.
 \end{aligned}
$$
Therefore
$$
 \int _{C_r} |h'(\xi)|^t  |d\xi | \succeq  d^{n(t-1)}\sum_{\xi\in p^{-{(n-M)}}(h(2))}|(p^{n-M})'(\xi)|^{-t},
$$
from which immediately follows that 
$$
\hat \b_{h}(t)\geq t-1+ \frac{\P(t)}{\log d}.
$$
On the other hand, $\hat \b_{h}(t)\leq \b_{h}(t)$ and, arguing exactly as before but with $C_r$ replaced by the full
circle $\{|\xi |=r\}$ and skipping the mixing argument based on the existence of the integer $M$ (which is no longer needed), it follows that
$$
\b_{h}(t)\le t-1+ \frac{\P(t)}{\log d}.
$$
\epf


\bibliographystyle{plain}


\end{document}